\documentclass[bj]{imsart}

\RequirePackage{amsthm,amsmath,amsfonts,amssymb,bbm}
\RequirePackage[numbers]{natbib}
\RequirePackage[colorlinks,citecolor=blue,urlcolor=blue]{hyperref}
\RequirePackage{graphicx}

\startlocaldefs
\newcommand{\ba}{\begin{eqnarray}}
\newcommand{\ea}{\end{eqnarray}}

\newcommand{\follmer}{{F\"ollmer}}
\newcommand{\cadlag}{{c\`adl\`ag }}

\theoremstyle{plain}

\newtheorem{theorem}{Theorem}[section]
\newtheorem{lemma}[theorem]{Lemma}
\newtheorem{corollary}[theorem]{Corollary}
\newtheorem{proposition}[theorem]{Proposition}
\theoremstyle{remark}
\newtheorem{remark}{Remark}
\newtheorem{definition}[theorem]{Definition}

\newtheorem*{example}{Example}

\endlocaldefs

\begin{document}

\begin{frontmatter}
\title{Quadratic variation and quadratic roughness}
\runtitle{Quadratic roughness}

\begin{aug}
\author[A]{\fnms{Rama} \snm{Cont}\ead[label=e1]{Rama.Cont@maths.ox.ac.uk}}
\and
\author[A]{\fnms{Purba} \snm{Das}\ead[label=e2,mark]{Purba.Das@maths.ox.ac.uk}}

\address[A]{Mathematical Institute,
University of Oxford.\\
\printead{e1,e2} }

\end{aug}

\begin{abstract}
We study the concept of quadratic variation of a continuous path along a sequence of partitions and its dependence with respect to the choice of the partition sequence. We introduce the concept of {\it quadratic roughness} of a path along a partition sequence and show that for H\"older-continuous paths satisfying this roughness condition, the quadratic variation along balanced partitions is invariant with respect to the choice of the partition sequence. Typical paths of Brownian motion are shown to satisfy this quadratic roughness property almost-surely along any partition with a required step size condition. Using these results we derive a formulation of the pathwise F\"ollmer-It\^o calculus which is invariant with respect to the partition sequence. We also derive an invarience of local time under quadratic roughness.
\end{abstract}

\begin{keyword} 
\kwd{Quadratic variation}
\kwd{Pathwise integration}
\kwd{Brownian motion} 
\kwd{It\^o calculus}  \kwd{Local time}
\kwd{Roughness} 
\end{keyword}

\end{frontmatter}

The concept of {\it quadratic variation} plays a central role in stochastic analysis and in the modern theory of stochastic integration \cite{dm,protter}. The quadratic variation of a (real-valued) random process $(X(t),t\in [0,T])$ with \cadlag \ sample paths is defined as the limit in the sense of (uniform) convergence in probability, of the sum of squared increments 
\ba\label{eq.qv}{ \sum_{\pi_n} \left(X(t^n_{k+1}\wedge t)-X(t^n_k\wedge t)\right)^t\left(X(t^n_{k+1}\wedge t)-X(t^n_k\wedge t)\right) }\ea
computed along a sequence of partitions $\pi^n=(0=t^{n}_0<t^n_1<\cdots<t^{n}_{N(\pi^n)}=T)$ with vanishing step size $|\pi^n|=\sup_{i=1,\cdots, N(\pi^n)}|t^n_i-t^n_{i-1}|\to 0$. The relevance of this notion, as opposed to $p-$variation, is underlined by the fact that large classes of random processes --such as Brownian motion and diffusion processes-- have finite quadratic variation, while at the same time possessing infinite $2-$variation.

Although quadratic variation for a stochastic process $X$ is usually defined as a limit in probability of \eqref{eq.qv}, it is essentially a pathwise property.
In his seminal paper {\it Calcul d'It\^o sans probabilit\'es} \cite{follmer1981}, Hans F\"ollmer introduced the class of c\`adl\`ag paths $X\in D([0,T],\mathbb{R})$ with finite quadratic variation along a sequence of partitions $\pi=(\pi^n)$, for which \eqref{eq.qv} has a limit with Lebesgue decomposition $[X]_\pi(t)=[X]^c(t)+\sum_{0\leq s\leq t}(\Delta X_s)^2$ and showed that
for $f\in C^2(\mathbb{R})$ one can define the integral $\int_0^. (\nabla f\circ X) d^\pi X$ as a pointwise limit of left Riemann sums along $(\pi^n)$:
\begin{equation}
 \int_0^T (\nabla f\circ X) d^\pi X :=\lim_{n\to\infty} \sum_{\pi^n}\nabla f(X(t)).(X(t^n_{i+1}\wedge T)-X(t^n_i\wedge T)),\label{eq.follmerintegral}\end{equation}
and
 this integral satisfies a change of variable formula:
 \begin{eqnarray}
 f( X(t)) = f( X(0)) + \int_0^t (\nabla f\circ X) d^\pi X + \frac{1}{2} \int_0^t \nabla^2 f ( X(s)). d [ X]^c_\pi
\nonumber \\+\sum_{s\in[0,t]}\Big( f( X(s))-f( X(s-))-\nabla f( X(s))\Delta X(s)\Big).\label{eq.ito}
 \end{eqnarray}
This `pathwise It\^o formula' may be used as a starting point for a purely pathwise construction of the It\^o calculus \cite{follmer1981,chiu2021} but, unlike the analogous theory for Riemann-Stieltjes or Young integrals, the construction in \cite{follmer1981} seems to depend on the choice of the sequence of partitions $(\pi^n)$: both the quadratic variation $[X]_\pi$ and the pathwise integral \eqref{eq.follmerintegral} are defined as limits along this sequence of partitions.
In fact, as shown by Freedman \cite[p. 47]{freedman}, 
for any continuous function $x$ one can construct a sequence of partitions $\pi$ such that $[x]_\pi=0$. This result was extended by Davis et al. \cite{davis2018} where they have  shown that 
given any continuous path $x$ and any increasing function $A$, one can construct a sequence of partitions $\pi$ such that $[x]_{\pi}=A$. These negative results seem to suggest that the dependence of $[x]_\pi$ on $\pi$ leaves no hope for the uniqueness of the quantities in Equation \eqref{eq.ito}.

On the other hand, as shown by L\'evy \cite{levy1940,levy1965} and Dudley \cite{dudley1973}, for typical paths of Brownian motion the sums in \eqref{eq.qv} converge to a unique limit along {\it any} sequence of partitions which are refining or whose mesh decreases to zero fast enough. Therefore there exists a large set of paths- containing all typical Brownian paths -for which one should be able to define the quantities in Equation \eqref{eq.ito} independent of the choice of the partition sequence $(\pi^n)_{n\geq 1}$ for a large class of sequences. 


We clarify these issues by investigating in detail the dependence of quadratic variation with respect to the sequence of partitions and deriving sufficient conditions for the stability of quadratic variation with respect to the choice of partition sequence. These conditions are related to an irregularity property of the path, which we call {\it quadratic roughness} (Def. \ref{def.rough}): this property requires cross-products of increments along the partition to average to zero at certain scales and is different from other notions of roughness such as H\"older roughness \cite{FrizHairer} or the concept of $\rho-$irregularity as put forth by Catellier and Gubinelli \cite{catellier2016}. H\"older roughness, like H\"older regularity, involves the amplitude of increments of a function, whereas our definition crucially involves the sign of the increments (or `phase' in the multidimensional case). The relation between quadratic roughness and $\rho-$irregularity \cite{catellier2016} is less clear. $\rho-$irregularity is based on the smoothness of the local time of a path; while our approach relies only on the existence of quadratic variation along certain  partition sequences and does not require the existence of a local time, it is possible that such properties would be implied by the existence of a smooth local time.

We present two main results in this paper.
First, we show that the quadratic roughness property is satisfied almost-surely by Brownian paths (Theorem \ref{thm.brownian}). 
Our second main result (Theorem \ref{main.theorem}) is an invariance result for pathwise quadratic variation: we show that for H\"older-continuous paths satisfying this irregularity condition, the quadratic variation along balanced partitions (Def. \ref{def.balance}) is invariant with respect to the choice of the partition sequences.
This leads to an {\it invariant} definition of quadratic variation, and a robust formulation of the pathwise F\"ollmer-It\^o calculus (Theorem \ref{thm.ito}) for such irregular paths.

Our results complement previous work on the pathwise approach to It\^o calculus
\cite{ananova2017,chiu2018,CF2010,davis2014,davis2014,davis2018,follmer1981,imkeller2015,kim2019,perkowski2015} by identifying a set of paths for which these results are robust to the choice of the sequence of partitions involved in the construction.
In contrast to the constructions in \cite{karandikar1983,karandikar1995,karandikar2014}, our approach does not rely on any probabilistic tools.

{\bf Outline.}
Section \ref{sec:definitions} recalls the definition of quadratic variation along a sequence of partitions, following \cite{chiu2018,follmer1981}.
Section \ref{sec:wellbalanced} defines the class of balanced partition sequences and discusses asymptotic comparability of such partitions.
Section \ref{sec:roughness} introduces the concept of {\it quadratic roughness} and explores some of its properties. In particular, we show that typical Brownian paths almost-surely satisfy this irregularity property (Theorem \ref{thm.brownian}).

Section \ref{sec:mainresult} shows that the quadratic roughness of a path is a sufficient condition for the invariance of quadratic variation with respect to the choice of partitions (Theorem \ref{main.theorem}). 
This result allows to give a definition of quadratic variation invariant with respect to the choice of the partition sequences (Proposition \ref{prop.qv}).
Section \ref{sec:ito} builds on these results to arrive at a robust formulation of the pathwise F\"ollmer-It\^o calculus.
Section \ref{sec:localtime} extends these results to pathwise local time.

\section{Quadratic variation along a sequence of partitions}\label{sec:definitions}

Let $T>0$. We denote
$D([0,T],\mathbb{R}^d)$ the space of $\mathbb{R}^d$-valued right-continuous functions with left limits (c\`adl\`ag functions),
 $C^{0}([0,T],\mathbb{R}^d)$ the subspace of continuous functions and,
for $0< \nu < 1,$ 
$C^{\nu}([0,T],\mathbb{R}^d)$ the space of H\"older continuous functions with exponent $\nu$:
$$ C^{\nu}([0,T],\mathbb{R}^d) =\left\{ x\in C^0([0,T],\mathbb{R}^d)\quad \Big|\; \sup_{(t,s)\in[0,T]^2, t\neq s }\frac{\|x(t)-x(s)\|}{|t-s|^{\nu}}<+\infty  \right\}\subset C^0([0,T],\mathbb{R}^d),$$
$${\rm and}\qquad C^{\nu-}([0,T],\mathbb{R}^d) =\mathop{\bigcap}_{0\leq \alpha< \nu}C^{\alpha}([0,T],\mathbb{R}^d). $$ 

We denote by $\Pi([0,T])$ the set of all finite partitions of $[0,T]$.
A sequence of { partitions} of $[0,T]$ is a sequence $(\pi^n)_{n\geq 1}$ of elements of $\Pi([0,T])$:
$$\pi^n=\left(0=t^{n}_0<t^n_1<\cdots<t^{n}_{N(\pi^n)}=T\right).$$ 
We denote $N(\pi^n)$ the number of intervals in the partition $\pi^n$ and
\ba |\pi^n|=\sup\{ |t^n_i-t^n_{i-1}|, i=1,\cdots,N(\pi^n)\},\quad \underline{\pi^n }=\inf\{ |t^n_i-t^n_{i-1}|, i=1, \cdots , N(\pi^n)\}, \ea 
the size of the largest (resp. the smallest) interval of $\pi^n$.

\begin{example} Let $k\geq 2$ be an integer. The $k$-adic partition sequence of $[0,T]$ is defined by $$ \pi^n=\bigg(t^n_j= \frac{j \ T}{k^n}, \qquad j=0,\cdots, k^n\bigg).$$We have $ \underline{\pi^n}=|\pi^n|= T/k^{n}.$\end{example}

\begin{example}[Lebesgue partition]Given $x\in D([0,T], \mathbb{R}^d)$ define $$ \lambda^n_0(x)=0,\quad{\rm and}\ \forall k\geq 1; \quad \lambda^n_{k+1}(x)=\inf\{ t\in ( \lambda^n_k(x),T],\quad \|x(t)-x( \lambda^n_k(x))\| \geq 2^{-n} \}$$and $N(\lambda^n(x))=\inf\{k\geq 1, \quad \lambda^n_k(x)=T \}.$We call the sequence $\lambda^n(x)=( \lambda^n_k(x) )$ the (dyadic) Lebesgue partition associated to $x$.\end{example}

\begin{definition}[Quadratic variation of a path along a sequence of partitions]\label{def.pathwiseQV}
Let $\pi^n=(0=t^{n}_0<t^n_1<\cdots<t^{n}_{N(\pi^n)}=T)$ be a sequence of partitions  of  $[0,T]$ with vanishing mesh
$|\pi^n|= \mathop{\sup}_{i=0,\cdots, N(\pi^n)-1} |t^n_{i+1}-t^n_i| \to 0$.
A c\`adl\`ag function $x \in D([0,T],\mathbb{R})$ is said to have finite
quadratic variation along the sequence of partitions
$(\pi^n)_{n\geq 1}$ if 
 the sequence of  measures
 \begin{equation*}
    \sum_{t^n_j \in \pi^n} ( x(t^n_{j+1}) - x(t^n_j))^2 \delta_{t^n_j}\quad 
   \end{equation*}
  converges weakly on $[0, T]$ to a limit measure $\mu$ such that
$t\mapsto [x]_\pi^c(t) = \mu([0,t]) - \sum_{0 < s \leq t} |\Delta x(s)|^2$
 is continuous and increasing.
The increasing function $[x]_\pi:[0,T]\to \mathbb{R}_+$ defined by \ba [x]_\pi(t)=\mu([0,t])=\lim_{n\to\infty} \sum_{\pi_n} (x(t^n_{k+1}\wedge t)-x(t^n_k\wedge t))^2  \label{eq:qv}\ea  is  called the {\it quadratic variation of }  $x$ along the sequence of partitions $\pi$. We denote $Q_\pi([0,T] ,\mathbb{R})$ the set of c\`adl\`ag paths with these  properties.
\end{definition}

$Q_\pi([0,T] ,\mathbb{R})$ is not a vector space (see e.g \cite{schied2016}). 
The extension of pathwise quadratic variation to vector-valued paths requires some care \cite{follmer1981}:
\begin{definition}[Pathwise quadratic variation for a vector valued path]\label{defn.qv.vector}
A c\`adl\`ag path  $x=(x^1,...,x^d)\in D([0,T],\mathbb{R}^d)$ is said to have finite quadratic variation along $\pi=(\pi^n)_{n\geq 1}$ if for all $i,j=1,\cdots,d$ we have
$x^i\in Q_\pi([0,T] ,\mathbb{R})$ and $x^i+x^j\in Q_\pi([0,T] ,\mathbb{R})$. We then denote $[x]_\pi\in D([0,T], S^+_d)$ the matrix-valued function defined by
$$[x]_\pi^{i,j}(t)=\frac{[x^i+x^j]_\pi(t)-[x^i]_\pi(t)-[x^j]_\pi(t)}{2}$$
where
$S^+_d$ is the set of symmetric semidefinite positive matrices.
We denote by $Q_\pi([0,T] ,\mathbb{R}^d)$ the set of functions satisfying these properties.
\end{definition}
For $x\in Q_\pi([0,T] ,\mathbb{R}^d)$, $[x]_\pi$ is a \cadlag function with values in $S^+_d$: $[x]_\pi\in D([0,T],S^+_d)$.

As shown in \cite{chiu2018}, the above definitions may be more simply expressed in terms of convergence of discrete approximations.
For continuous paths, we have the following characterization \cite{cont2012,chiu2018} for quadratic variation:
\begin{proposition}\label{prop.cont.qv}
$x\in C^0([0,T],\mathbb{R}^d)$ has finite quadratic variation along partition sequence $\pi=(\pi^n,n\geq 1)$ if and only if the sequence of functions $\left([x]_{\pi^n}, \; n\geq 1\right)$ defined by
$$[x]_{\pi^n}(t):=\sum_{t^n_j\in \pi^n}\left(x(t^n_{j+1}\wedge t)-x(t^n_{j}\wedge t)\right)^t \left(x(t^n_{j+1}\wedge t)-x(t^n_{j}\wedge t)\right),$$
converges uniformly on $[0,T]$ to a continuous (increasing) function $[x]_{\pi} \in C^0([0,T],S^+_d)$.
\end{proposition} 
 
The notion of quadratic variation along a sequence of partitions is different from the p-variation for $p=2$. The p-variation involves taking a supremum  over {\it all} partitions, whereas quadratic variation is a limit taken along a specific partition sequence $(\pi^n)_{n\geq 1}$. In general $[x]_\pi$ given by \eqref{eq:qv} is smaller than the p-variation for $p=2$. In fact, for diffusion processes, the typical situation is that p-variation is (almost-surely) infinite for $p=2$ \cite{dudley2011,taylor1972} while the quadratic variation is finite for sequences satisfying some mesh size condition. 
For instance, typical paths of Brownian motion have finite quadratic variation along any sequence of partitions with mesh size $o(1/\log n)$ \cite{dudley1973,delavega1974} while simultaneously having infinite p-variation almost surely for $p\leq 2$ 
\cite[p. 190]{levy1965}:
$$ \inf_{\pi\in \Pi(0,T)} \sum_{\pi} |W(t_{k+1})- W(t_k)|^2 = 0,\qquad {\rm while}\qquad\sup_{\pi\in \Pi(0,T)} \sum_{\pi} |W(t_{k+1})- W(t_k)|^2 = \infty$$
almost-surely.

The quadratic variation of a path along a sequence of partitions strongly depends on the chosen sequence. In fact, as shown by Freedman \cite[p. 47]{freedman}, given any continuous functions, one can always construct a sequence of partitions along which the quadratic variation is zero. This result was extended by Davis et al. \cite{davis2018} who show that, given any continuous path $x\in C^0([0,T],\mathbb{R})$ and any increasing function $A:[0,T]\to\mathbb{R}_+$ one can construct a partition sequence $\pi$ such that $ [x]_{\pi}=A$.
Notwithstanding these negative results, we shall identify a class of paths $x$ for which $[x]_\pi$ is uniquely defined across the class of {\it balanced} partition sequences, which we now define.

\section{Balanced partition sequences}\label{sec:wellbalanced}
One difficulty in comparing quadratic variation along two different partition sequences is the lack of uniform bounds on the partition intervals and the lack of comparability between two partitions. 
In this section, we introduce the class of {\it balanced} partition sequences which allow such bounds.

We shall say two (real) sequences $a=(a_n)_{n\geq 1}$ and $b=(b_n)_{n\geq 1}$ are asymptotically comparable, denoted $a_n\asymp b_n$, if $|a_n|= O(|b_n|)$ and $|b_n|=O(|a_n|)$. If both sequences are strictly positive then
$$a_n\asymp b_n\quad \iff\quad\limsup_{n\to\infty} \frac{|b_n|}{|a_n|}<\infty\quad{\rm and}\quad \limsup_{n\to\infty} \frac{|a_n|}{|b_n|}<\infty.$$
\[\iff \exists M_0<\infty\; \text{ such that: } \; \forall n\in \mathbb{N}: \quad \frac{|b_n|}{|a_n|}<M_0\quad{\rm and}\quad  \frac{|a_n|}{|b_n|}<M_0.\]

\subsection{Definition and properties}
\begin{definition}[Balanced partition sequence]\label{def.balance} 
Let $\pi^n=\left(0=t^{n}_0<t^n_1<\cdots<t^{n}_{N(\pi^n)}=T\right)$ be a sequence of partitions of interval $[0,T]$ and $$\underline{\pi^n}=\inf_{i=0,\cdots, N(\pi^n)-1} |t^n_{i+1}-t^n_i|,\qquad |\pi^n|=\sup_{i=0,\cdots, N(\pi^n)-1} |t^n_{i+1}-t^n_i|. $$
We say $(\pi^n)_{n\geq 1}$ is {\em balanced partition sequence} if 
\begin{equation}\exists\; c>0,\;  \forall n\geq 1, \quad \frac {|\pi^n|}{\underline{\pi^n}}\leq c.\label{eq.balance}\end{equation}
\end{definition}
This condition means that all intervals in the partition sequence $\pi^n$ are asymptotically comparable. Note that, since $\underline{\pi^n}N(\pi^n)\leq T$, any balanced sequence of partitions satisfies the following inequality.
\begin{equation}
\underline{\pi^n}\leq |\pi^n|\leq c \ \underline{\pi^n} \leq \frac{cT}{N(\pi^n)}.\label{eq.wellbalanced}
\end{equation}
We denote by $\mathbb{B}([0,T])$, the set of all balanced partition sequences of $[0,T]$. 
\begin{proposition}[Properties of balanced partition sequence]\label{prop.well.balanced}
Let $\pi=(\pi^n)_{n\geq 1}$ be a sequence of partitions of $[0,T]$ with mesh $|\pi^n|\to 0$. Then:
\begin{enumerate}
 \item[(i)] $\pi\in \mathbb{B}([0,T])\iff \liminf_{n\rightarrow \infty} N(\pi^n)\underline{\pi^n} > 0$ and $\limsup_{n\rightarrow \infty} N(\pi^n)|\pi^n| < \infty. $
 \item[(ii)] Let $N(\pi^n,t_1,t_2)$ be the number of partition points of $\pi^n$ in $[t_1,t_2]$. If $\pi\in \mathbb{B}([0,T])$ then for any $h>0$,
$$\limsup_{n\rightarrow\infty} \frac{\sup_{t\in[0,T-h]} N(\pi^n,t,t+h)}{\inf_{t\in[0,T-h]} N(\pi^n,t,t+h)}<\infty. $$
\item[(iii)] If $\pi=\left(\pi^n;\;n\geq 1\right) \in \mathbb{B}([0,T])$ then
\begin{equation}
 \limsup_n \frac{N(\pi^{n+1})}{N(\pi^n)}<\infty \iff \limsup_n \frac{|\pi^n|}{|\pi^{n+1}|}<\infty \iff \limsup_n \frac{\underline{\pi^n}}{\underline{\pi^{n+1}}}<\infty.\label{eq.growth}
\end{equation} 
\item[(iv)] If $g\in C^1([0,T],\mathbb{R})$ is strictly increasing with $\inf g' > 0$ then the image under $g$ of a balanced partition sequence of $[0,T]$ is also a balanced partition sequence of $\left[g(0),g(T)\right].$
\end{enumerate}
\end{proposition}
The proof of this Proposition is given in Appendix \ref{appendix.2.2}.

\begin{definition}[Asymptotic comparability of balanced partitions]\label{def.comparability} 
We will say that two balanced partition sequences $\tau=(\tau^n)_{n\geq 1}$ and $\sigma=(\sigma^n)_{n\geq 1}$ are (asymptotically) comparable if 
\begin{equation}0< \mathop{\liminf}_{n\to\infty}\frac{|\sigma^n|}{|\tau^n|}\leq\mathop{\limsup}_{n\to\infty}\frac{|\sigma^n|}{|\tau^n|}<\infty\label{eq.stepcomparability}.\end{equation}
Since the partition sequences are balanced (not true for a general partition sequence), an equivalent condition will be:
\begin{equation}
0< \mathop{\liminf}_{n\to\infty}\frac{N(\sigma^n)}{N(\tau^n)}\leq\mathop{\limsup}_{n\to\infty}\frac{N(\sigma^n)}{N(\tau^n)}<\infty\label{eq.Ncomparability}.\end{equation}
We denote $\tau \asymp \sigma$ (or $\tau^n \asymp \sigma^n$).
\end{definition}
Note that for general (not balanced) sequences of partitions Inequality \eqref{eq.stepcomparability} neither implies nor  is implied by Inequality \eqref{eq.Ncomparability}: this is purely a  consequence of \eqref{eq.wellbalanced}. If $\tau \asymp \sigma$ then the number of partition points of $\tau^n$ (respectively $\sigma^n$) in any consecutive partition points of $\sigma^n$ (respectively $\tau^n$) remains bounded as $n\to\infty$.

The following Lemma shows how one can adjust the rate at which the mesh of a balanced sequence converges to zero.
\begin{lemma}\label{subseq.asymp.comp}
Let $\tau=(\tau^n)_{n\geq 1}$ and $\sigma=(\sigma^n)_{n\geq 1}$ be two balanced partition sequences of $[0,T]$ with $\limsup_n \frac{|\sigma^n|}{|\tau^n|} < 1$ and mesh $|\tau^n|\xrightarrow[]{n\rightarrow \infty} 0$. 
\begin{enumerate}
 \item[(i)]
 There exists a subsequence $(\tau^{k(n)})_{n\geq 1}$ of $\tau$ such that: 
$$|\tau^{k(n)}|\to 0 \quad \text{ and, } \quad \limsup_n \frac{| \sigma^{n}|}{|\tau^{k(n)}|} \geq 1.$$
\item[(ii)] Furthermore if we also assume $\limsup_n\frac{|\tau^n|}{| \tau^{n+1}|} < \infty$, then there exists a subsequence $(\tau^{k(n)})_{n\geq 1}$ of $\tau$ which is asymptotically comparable to $\sigma$: $\tau^{k(n)}\asymp \sigma^n$.
\item[(iii)] There exists $r:\mathbb{N}\mapsto \mathbb{N}$ with $\lim_{n\to \infty}r(n)=\infty$ such that $$\limsup_n \frac{| \sigma^{r(n)}|}{|\tau^{n}|} \geq 1.$$
\end{enumerate}
\end{lemma}
Note that $r:\mathbb{N}\mapsto \mathbb{N}$ in Lemma \ref{subseq.asymp.comp} (iii) is not necessarily  injective i.e. $(\sigma^{r(n)},n\geq 1)$ is not necessarily  a subsequence of $(\sigma^{n},n\geq 1)$.

\begin{proof}
Denote the partition points of $\tau^n$ and $\sigma^n$ respectively by $\big(t_{k}^n; \; k=0,\cdots, N(\tau^n)\big)$ and $\big(s_{l}^n; \; l=0,\cdots, N(\sigma^n)\big)$.
\\\textit{Proof of (i)}: From the assumption, $\limsup_n \frac{|\sigma^n|}{|\tau^n|} < 1$ which implies $\liminf_n\frac{|\tau^n|}{|\sigma^n|} > 1$. Then there exists $N_0\in \mathbb{N}$ such that for all $n\geq N_0,\; \frac{|\tau^n|}{|\sigma^n|} > 1$. Since we are only concerned about the limiting behaviour as $n\rightarrow \infty$ we will only consider $n>N_0$.
We define $k(n)$ as follows.
\begin{equation}
 \label{eq.kn}
k(n) = \inf\{ k\geq n; \quad |\tau^k| \leq |\sigma^n| \} <\infty\quad {\rm since,}\quad |\tau^k|\mathop{\to}^{k\to\infty} 0.\end{equation}
We now consider the subsequence $(\tau^{k(n)})_{n\geq 1}$ of 
$\tau$. From the definition of $k(n)$:
$$\limsup_n \frac{| \sigma^{n}|}{|\tau^{k(n)}|} \geq 1.$$
\textit{Proof of (ii)}: If $\limsup_n\frac{|\tau^n|}{|\sigma^n|} <\infty$ we set $k(n)=n$; otherwise if $\limsup_n\frac{|\tau^n|}{|\sigma^n|} =+\infty$,
define $k(n)$ as in Equation \eqref{eq.kn}. Now for $i= 1,\cdots, N(\sigma^n)$ define $j(i,n)$ as follows.
$$ j(i,n) = \inf\{ j\geq 1, \quad t^{k(n)}_{j}\in (s^n_{i},s^n_{i+1}] \}.$$ 
Then we have: $$ t^{k(n)}_{j(i,n)-1} \leq s^n_k < t^{k(n)}_{j(i,n)} < \cdots < t^{k(n)}_{j(i+1,n)-1} \leq s^n_{i+1} < t^{k(n)}_{j(i+1,n)} .$$
If for some $i,$ \; $|j(i+1,n)- j(i,n)| \rightarrow \infty$ as $n\to\infty$ then, from the above construction of $k(n)$ and using the well balanced property of $\sigma^n $ and $\tau^{k(n)}$ we have: $\limsup_n \frac{|\sigma^n|}{|\tau^{k(n)}|}\rightarrow \infty $ and $\limsup_n \frac{|\sigma^n|}{|\tau^{k(n)-1}|}<1$. Hence, $\limsup_n \left(\frac{|\sigma^n|}{|\tau^{k(n)}|}-\frac{|\sigma^n|}{|\tau^{k(n)-1}|}\right) = \limsup_n \frac{|\sigma^n|}{|\tau^{k(n)-1}|}\left[\frac{|\tau^{k(n)-1}|}{|\tau^{k(n)}|} -1 \right]$ $\rightarrow \infty $ which is a contradiction because of our assumption. Hence the cluster size $j(i+1,n)- j(i,n)$ is uniformly bounded:
$$\forall i,n\geq 1, \;\; \exists M, \quad \text{such that, }\; |j(i+1,n)- j(i,n)|\leq M< \infty.$$
So there exists a constant $c_0$ such that
\begin{equation}
 1\leq \liminf_{n} \frac{|\sigma^n|}{|\tau^{k(n)}|} \leq \limsup_{n} \frac{|\sigma^n|}{|\tau^{k(n)}|}\leq c_0 < \infty . 
\end{equation}
Therefore $(\tau^{k(n)})_{n\geq 1}$ and $(\sigma^{n})_{n\geq 1}$ are (asymptotically) comparable.
\\ \textit{Proof of (iii)}: Since $\limsup_n \frac{|\sigma^n|}{|\tau^n|} < 1$, the set $\{ n\geq 1,\ \frac{|\sigma^n|}{|\tau^n|} \geq 1\}$ is finite and the set,
$$ A=\{ n\geq 1,\ \frac{|\sigma^n|}{|\tau^n|} < 1\}$$ is infinite.
Now from the assumption there exists $N_0\in \mathbb{N}$ such that for all $n\geq N_0,\; \frac{|\tau^n|}{|\sigma^n|} > 1$. 
Now  for $n\leq N_0$ set $r(n)=1$, for $n>N_0$ and  $n\notin A$,  set $r(n)=n$ and 
$${\rm for}\quad n>N_0;\; n\in A \qquad r(n)=\sup\{r\leq n,\; |\sigma^r|>|\tau^n| \}<\infty. $$
Then,
$$ r(n)\leq n \quad and, \quad \mathop{\limsup}_{n\to\infty} \frac{|\sigma^{r(n)}|}{|\tau^n|} = \mathop{\limsup}_{n\in A} \frac{|\sigma^{r(n)}|}{|\tau^n|} \geq 1.$$
\end{proof}

\subsection{Quadratic variation along balanced partition sequences}
If a path has quadratic variation along a sequence of partitions, then it also has (the same) quadratic variation along any sub-sequence. This simple remark has interesting implications when the partition sequences are balanced: comparing the sum of squared increments along the original sequence with the sum along a sub-sequence (with finer mesh) we obtain that, under some scaling conditions on the mesh, cross-products of increments along the finer partition average to zero across the coarser partition.
\begin{lemma}[Averaging property of cross-products of increments]\label{lemma.averaging}
Let $x\in  C^{\alpha}([0,T],\mathbb{R}^d)$ for some $\alpha>0$ and $\sigma^n=\left(0=s^n_0<s^n_1<\cdots< s^n_{N(\sigma^n)}=T\right)$ be a balanced sequence of partitions of $[0,T]$ such that $x\in Q_\sigma([0,T],\mathbb{R}^d)$. 
Let $\kappa> \frac{1}{2\alpha}$ and $(\sigma^{l_n})_{n\geq 1}$ a subsequence of $\sigma^n$ with $|\sigma^{l_n}|= O\left(|\sigma^n|^{\kappa}\right)$.
For $k= 1 ,\cdots, N(\sigma^n)$ define $p(k,n) = \inf\{ m\geq 1: \quad s^{l_n}_{m}\in (s^n_{k},s^n_{k+1}] \}$ . Then 
$$\sum_{k=1}^{N(\sigma^n)} \sum_{p(k,n)\leq i\neq j < p(k+1,n)-1} \bigg(x(s^{l_n}_{i+1})-x(s^{l_n}_{i})\bigg)^t \bigg(x(s^{l_n}_{j+1})-x(s^{l_n}_{j})\bigg) \xrightarrow[]{n\to \infty} 0.$$
\end{lemma}
\begin{proof}
We provide the proof for $d=1$. The extension to $d>1$ is straightforward extension of $1$-dimensional case.
Let $\sigma^{l_n}$ be a sub-sequence of $\sigma^n$ satisfying $|\sigma^{l_n}|= O\left(|\sigma^n|^{\kappa}\right)$. Denote, \[[x]_{\sigma^{n}}(t) = \sum_{k=1}^{N(\sigma^n)-1} \left(x(s_{k+1}^n\wedge t)-x(s_{k}^n\wedge t)\right)^t\left(x(s_{k+1}^n\wedge t)-x(s_{k}^n\wedge t)\right),\quad \text{and,}\]
\[[x]_{\sigma^{l_n}}(t) =\sum_{s^{l_n}_k\in \sigma^{l_n}} \left(x(s_{k+1}^{l_n}\wedge t)-x(s_{k}^{l_n}\wedge t)\right)^t\left(x(s_{k+1}^{l_n}\wedge t)-x(s_{k}^{l_n}\wedge t)\right). \]
Then $\bigg|[x]_{\sigma^{l_n}}(t)-[x]_{\sigma^{n}}(t)\bigg|\rightarrow 0.$
Grouping the points of $\sigma^{l_n}$ along partition points of $\sigma^n$, we obtain:
$$\bigg|[x]_{\sigma^{n}}(T)-[x]_{\sigma^{l_n}}(T)\bigg|=\left|\sum_{\sigma^n}\left(x(s^n_{i+1})-x(s^n_{i})\right)^2 - \sum_{\sigma^{l_n}}\left(x(s^{l_n}_{i+1})-x(s^{l_n}_{i})\right)^2\right|$$
$$=\left|\sum_{\sigma^n}\bigg(\left(x(s^n_{i+1})-x(s^n_{i})\right)^2 - \sum_{j=p(i,n)}^{p(i+1,n)-2}\left(x(s^{l_n}_{j+1})-x(s^{l_n}_{j})\right)^2\bigg)+ \sum_{k=1}^{N(\sigma^n)} \left(x(s^{l_n}_{p(k,n)})-x(s^{l_n}_{p(k,n)-1})\right)^2\right|$$ 
$$\geq \left|\sum_{\sigma^n}\bigg( \left(x(s^n_{i+1})-x(s^n_{i})\right)^2 - \sum_{j=p(i,n)}^{p(i+1,n)-2}\left(x(s^{l_n}_{j+1})-x(s^{l_n}_{j})\right)^2\bigg) \right|- \sum_{k=1}^{N(\sigma^n)} \left(x(s^{l_n}_{p(k,n)})-x(s^{l_n}_{p(k,n)-1})\right)^2.$$ 
Using the $\alpha-$H\"older continuity of   $x$, the last term in the above equation is bounded above by $\sum_{k=1}^{N(\sigma^n)} C |\sigma^{l_n}|^{2\alpha} \leq CN(\sigma^n) |\sigma^{l_n}|^{2\alpha}$. Now using the balanced property of $\sigma^{l_n}$ (subsequence of a balanced sequence of partitions is also balanced), we further get the above bounded as: $$ \sum_{i=1}^{N(\sigma^n)} C|\sigma^{l_n}|^{2\alpha} \leq \frac{C_1 N(\sigma^n)}{N(\sigma^{l_n})^{2\alpha}} \leq C_2 \times N(\sigma^n)^{1-2\alpha \kappa} \mathop{\to}^{n\to\infty} 0,$$ since $1-2\kappa \alpha <0$. 
So writing the first term of the previous equation explicitly we finally obtain,
$$\lim_{n\rightarrow \infty}\left| \sum_{k=1}^{N(\sigma^n)} \sum_{p(k,n)\leq i\neq j< p(k+1,n)-2} \left(x(s^{l_n}_{i+1})-x(s^{l_n}_{i})\right) \left(x(s^{l_n}_{j+1})-x(s^{l_n}_{j})\right)\right|$$ $$\leq \lim_{n\rightarrow \infty} \left| [x]_{\sigma^{n}}-[x]_{\sigma^{l_n}}\right|= 0.$$
\end{proof}

\section{Quadratic roughness}\label{sec:roughness}

\subsection{Quadratic roughness along a sequence of partitions}

Lemma \ref{lemma.averaging} shows that if a function has finite quadratic variation along a balanced partition sequence, then the cross-products of the increments along any subsequence with {\it sufficiently small mesh} average to zero  along the original (coarser) sequence. 
Intuitively, this means that there is enough cancellation across neighbouring increments such that their cross-products average to zero under coarse-graining. This can only occur if the increments over any scale have alternating signs, which is an indicator of the `roughness' of the function itself.
We will now introduce a slightly extended version of this property, which we call {\it quadratic roughness}, and show that this property plays a crucial role in the stability of quadratic variation with respect to the choice of partition.

\begin{definition}[Super-sequence]\label{superseq}
We call $d^n=\pi^{r(n)}$  a {\it super-sequence} of $\pi= (\pi^n)_{n\geq 1}$ if the map $r: \mathbb{N} \to \mathbb{N}$ is non-decreasing and $k\geq r(k)$ for all $k \in \mathbb{N}$.

\end{definition}

\begin{definition}[Quadratic roughness]\label{def.rough} Let $\mathbb{T} = (\mathbb{T}^n)_{n\geq 1}$ be the dyadic partition of $[0,T]$ and
$\pi^n=\left(0=s^{n}_0<s^n_1<\cdots <s^{n}_{N(\pi^n)}=T\right)$ be a balanced sequence of partitions of $[0,T]$ with vanishing mesh  $|\pi^n|\to 0$.  
We say that $x\in C^0([0,T],\mathbb{R}^d)\cap Q_{\mathbb{T}}([0,T],\mathbb{R}^d)$ has the {\rm quadratic roughness} property with coarsening index $0< \beta<1$ along $\pi$ on $[0,T]$ if there exists a  subsequence or super-sequence 
$d^n= \left(0=t^n_1<t^n_2<\cdots<t^n_{N(d^n)}=T\right)$ of $\mathbb{T}$ with the following properties:
\begin{itemize}
    \item [(i)] $|d^n|^{\beta}=O\left(|\pi^n|\right)$ and,
    \item [(ii)] for all $t\in [0,T]$ : \[\sum_{j=1}^{N(\pi^n)-1}\sum_{t^n_i\neq t^n_{i'}\in (s^n_j,s^n_{j+1}]} \bigg(x(t^n_{i+1}\wedge t)-x(t^n_{i}\wedge t)\bigg)^t\bigg(x(t^n_{i'+1}\wedge t)-x(t^n_{i'}\wedge t)\bigg)\xrightarrow[]{n\to\infty} 0.\nonumber\label{eq.P1}\]
\end{itemize}
We denote by $R^\beta_\pi([0,T],\mathbb{R}^d)$ the set of paths satisfying this quadratic roughness property.
\end{definition}
In other words, the quadratic roughness property states that cross-products of increments along the dyadic partition $d^n$  average to zero when grouped along $\pi^n$.
Note that, since $\beta<1$, the number of terms in the inner sum in (ii) grows to infinity as $n$ grows, so (ii) is the result of compensation across terms, reminiscent of the law of large numbers.

\begin{remark}[Choice of reference partition]
In the quadratic roughness definition (\ref{def.rough}), we have used the dyadic partition as a `reference partition' to which other (balanced) partitions are compared. In fact, as will become clear in the proofs below, the dyadic partition may be replaced by any other balanced sequence of partitions $\sigma$ with vanishing mesh  $|\sigma^n|\to 0$ satisfying  $\sup_n\frac{|\sigma^n|}{|\sigma^{n+1}|}<\infty$ without changing the statements of any of the results.
\end{remark}

As a consequence of the quadratic roughness (Definition \ref{def.rough}) if $x\in R^\beta_\pi([0,T],\mathbb{R}^d)$ for some sequence of partitions $\pi$ of $[0,T]$, then $x\in Q_{\mathbb{T}}([0,T],\mathbb{R}^d)$ (but not necessarily $x\in Q_{\pi}([0,T],\mathbb{R}^d)$).

\begin{proposition}[Properties of quadratic roughness]\label{prop.rough.fn} Let $\pi= (\pi^n)_{n\geq 1}$ be a balanced partition sequence of $[0,T]$ with vanishing mesh (i.e. $|\pi^n| \to 0$) and
 $x\in R^\beta_\pi([0,T],\mathbb{R}^d)$  with $0 < \beta<1$. Then:
 \begin{enumerate}
\item  For any interval $I\subset [0,T]$, the path $x$ also has the quadratic roughness property on $I$ along stopped partition $\pi_I = (\pi^n_I)_{n\geq 1}= (\pi^n \cap I)_{n\geq 1}$. ie. $x\in R^\beta_{\pi_I}(I,\mathbb{R}^d)$. 
\item  For any subsequence/super-sequence $\tau^n=\pi^{k(n)}$ of $\pi$, we have $x\in R^\beta_\tau([0,T],\mathbb{R}^d)$.
\item  For any  $\lambda\in \mathbb{R}$, $\lambda x \in R^\beta_\pi([0,T],\mathbb{R}^d)$.
\item If $y$ is a function of finite variation then, $x+y\in R^\beta_\pi([0,T],\mathbb{R}^d)$.
\par Furthermore, if $d=1$ and if $y$ is a function with finite $p$-variation with $p<2$  then, $x+y\in R^\beta_\pi([0,T],\mathbb{R})$.
\item For  $\gamma\in[\beta,1)$, $R^\beta_\pi([0,T],\mathbb{R}^d)\subset R^\gamma_\pi([0,T],\mathbb{R}^d)$.
\end{enumerate}
\end{proposition}
The proof of this Proposition is given in Appendix \ref{appendix.A2}.

\subsection{Quadratic roughness of Brownian paths}\label{sec.Brownianrough}
We will now show that the quadratic roughness property is satisfied almost-surely by typical sample paths of Brownian motion.

\begin{theorem}[Quadratic roughness of Brownian paths]\label{thm.brownian}
Let $W$ be a Wiener process on a probability space $(\Omega, {\cal F},\mathbb{P})$, $T>0$ and
$(\pi^n)_{n\geq 1}$  a balanced sequence of partitions of $[0,T]$ with \begin{equation}
(\log n)^2 |\pi^n| \mathop{\to}^{n\to\infty}0. \label{eq.logscaling}\end{equation}
Then  the sample paths of $W$ almost-surely satisfy the quadratic roughness property for any $0< \beta < 1$:
$$\forall \beta\in (0,1),\quad \mathbb{P}\left( \ W\in R^\beta_\pi([0,T],\mathbb{R})\ \right)=1.$$
\end{theorem}
\begin{proof}
Let $W$ be a Wiener process on a probability space $(\Omega, {\cal F},\mathbb{P})$, which we take to be the canonical Wiener space without loss of generality i.e $\Omega=C^0([0,T],\mathbb{R}), W(t,\omega)=\omega(t)$.

Take $\beta\in (0,1)$. We know $\pi^n=(0=t^{n}_0<t^n_1<\cdots<t^{n}_{N(\pi^n)}=T)$ be a balanced sequence of partitions of $[0,T]$ satisfying Equation \eqref{eq.logscaling}. Now we will define a subsequence $\left(\mathbb{T}^{l_n}\right)$ of $\left(\mathbb{T}^n \right)$ such that $ |\mathbb{T}^{l_n}|^\beta =O(|\pi^n|)$.

\par Now, if $|\mathbb{T}^n|^\beta = O(|\pi^n|)$, then take $l_n=n$. Otherwise if $\limsup_n  \frac{|\mathbb{T}^n|^\beta }{|\pi^n|}=\infty$, then we define $l_n$ as follows,
\[l_n=\inf \{l\geq n: |\pi^n|\geq |\mathbb{T}^{l}|^{\beta}\}<\infty \qquad \text{since } |\mathbb{T}^l|^\beta=\frac{1}{2^{l\beta}}\xrightarrow[]{l\to \infty} 0.\label{eq.kkn1} \]
So from the construction of $l_n$ we then get:
\[|\mathbb{T}^{l_n}|\leq |\pi^n|^{1/\beta} < |\mathbb{T}^{l_n-1}|.\] Since the subsequence $(\mathbb{T}^{l_n})$ is also balanced, there exists constants $c_1$ and $c_2$ such that $c_1N(\mathbb{T}^{l_n})\geq N(\pi^{n})^{1/\beta} > c_2N(\mathbb{T}^{l_n}).$ 
Hence the dyadic subsequence $\mathbb{T}^{l_n}$ satisfies
\[|\mathbb{T}^{l_n}|^\beta  = O(|\pi^n|)\]

\par Now we will show that $\mathbb{T}^{l_n}$ satisfies Condition (ii) of Definition \ref{def.rough}. Define\\  $a^n_{ii'}= \sqrt{(t^{l_n}_{i+1}-t^{l_n}_{i})(t^{l_n}_{i'+1}-t^{l_n}_{i'})}$ if  $\exists \; j\in \{1,2,\cdots ,N(\pi^n)\}$ such that $p(n,j-1)\leq i\neq i'< p(n,j)$. Otherwise, set $a^n_{ii'}=0$. 
Let \[S_\pi(\mathbb{T}^{l_n},W)=\sum_{j=1}^{N(\pi^n)}\sum_{p(n,j-1)\leq i\neq i'< p(n,j)} \left(W(t^{l_n}_{i+1})-W(t^{l_n}_{i})\right)^t\left(W(t^{l_n}_{i'+1})-W(t^{l_n}_{i'})\right)\]\[=\sum_{i,i'=1}^{N(\mathbb{T}^{l_n})} a^n_{ii'} X^n_iX^n_{i'},\]
\[\text{where,}\quad X^n_i = \frac{W(t^{l_n}_{i+1})-W(t^{l_n}_{i}) }{\sqrt{t^{l_n}_{i+1}-t^{l_n}_{i}}} \sim N(0,1) \quad\text{ are IID variables for }\;  i=0,\cdots, N(\mathbb{T}^{l_n})-1.\]
Now let, \[\Lambda^2 =\sum_{1\leq i, i'\leq N(\mathbb{T}^{l_n})} (a^n_{ii'})^2 = \sum_{j=1}^{N(\pi^n)}\sum_{p(n,j-1)\leq i\neq i'< p(n,j)}(a_{ii'}^n)^2 \]
\[\leq \sum_{j=1}^{N(\pi^n)}\sum_{p(n,j-1)\leq i\neq i'< p(n,j)} \Delta t_i^{l_n}  \Delta t_{i'}^{l_n}\leq \sum_{j=1}^{N(\pi^n)} |\pi^n|^2 \leq |\pi^n|\sum_{j=1}^{N(\pi^n)} |\pi^n|\leq cT|\pi^n|.\]
The last inequality is due to the fact that $\pi$ is a balanced sequence. The Hanson-Wright inequality \cite{hanson1971} then implies that there exists constants $C_1$ and $C_2$ such that
$$\forall\delta>0,\quad \forall n\geq 1,\qquad \mathbb{P}\bigg(\left|S_\pi(\mathbb{T}^{l_n},W)\right|>\delta \bigg) \leq 2\exp(-\min\{C_1 \frac{\delta}{\sqrt{ |\pi^n|}}, C_2 \frac{\delta^2}{|\pi^n|} \}).$$
Since $|\pi^n|(\log n)^2\rightarrow 0 $   for large $n$, the upper bound is determined by the first term 
$\exp\left(-C_1 \frac{\delta}{\sqrt{|\pi^n|}}\right).$
If we denote   $\varepsilon_n^2= |\pi^n|(\log n)^2$ then $\varepsilon_n\rightarrow 0$ and we can rewrite this bound as
\begin{equation}
\mathbb{P}\bigg(\left|S_\pi(\mathbb{T}^{l_n},W)\right|>\delta \bigg) 
\leq 2\exp(-\min\{ \frac{C_1 \delta \ \log n}{\varepsilon_n},  \frac{C_2\delta^2\  (\log n)^2}{\varepsilon_n^2} \})\ \mathop{\leq} \frac{2 C}{n^{C_1\delta/\varepsilon_n}}. \label{eq.HR} \end{equation}

The series $\sum_n \frac{1}{n^{C_1\delta/\varepsilon_n}} < \infty$ is absolutely convergent. 
So we can apply the Borel-Cantelli lemma to obtain for each $\delta>0$ a set $\Omega_\delta$ with $\mathbb{P}(\Omega_\delta)=1$ and $N_\delta\in \mathbb{N}$ such that
$$\forall \omega\in \Omega_\delta,\quad \forall n\geq N_\delta, \quad \left|S_\pi(\mathbb{T}^{l_n},W)(\omega)\right|\leq \delta.$$
Now if we set $$\Omega_\pi=\Omega_0 \cap\bigg(\mathop{ \cap}_{m\geq 1}\Omega_{1/m} \bigg) \qquad {\rm then}\quad \mathbb{P}(\Omega_\pi)=1$$
and
for paths in $\Omega_\pi$ we have $S_\pi(\mathbb{T}^{l_n},\omega)\to 0$ simultaneously :
$$\forall \omega\in \Omega_\pi,\qquad S_\pi(\mathbb{T}^{l_n},\omega)\mathop{\to}^{n\to\infty} 0,$$
Therefore, $\mathbb{P}\left( \ W\in R^\beta_\pi([0,T],\mathbb{R})\ \right)=1.$
\end{proof}
The proof above uses independence of increments, which then implies that the cross-products of increments averages to zero due to a concentration inequality. However, the
quadratic roughness property is a pathwise property, and also holds for  classes  of stochastic processes with dependent increments, as the following example shows:
\begin{example}[Quadratic roughness of mixed Brownian motion]
Let $H>\frac{1}{2}$ and $\delta >0$ and 
$$M^{H,\delta}=  B+\delta B^H$$  where
$B$ is a Brownian motion and $B^H$ is a fractional Brownian motion with Hurst parameter $H$ on a probability space $(\Omega, {\cal F},\mathbb{P})$. Then for any
  balanced sequence $(\pi^n)_{n\geq 1}$ of partitions of $[0,T]$ with 
\begin{equation}
(\log n)^2|\pi^n|\ \mathop{\to}^{n\to\infty}0.
\end{equation}
  the sample paths of $M^{H,\delta}$ almost-surely satisfy the quadratic roughness property on $[0,T]$: 
$$\forall \beta\in (0,1),\quad \mathbb{P}\left( \ M^{H,\delta}\in R^\beta_\pi([0,T],\mathbb{R})\ \right)=1.$$
This is an application of Proposition \ref{prop.rough.fn} (iv) and Theorem \ref{thm.brownian}.
\end{example}

\section{Uniqueness of quadratic variation along balanced partitions}\label{sec:mainresult}

\subsection{Main result}
The following lemma shows that the quadratic roughness property is a necessary condition for the stability of quadratic variation with respect to the choice of the partition sequences:
\begin{lemma}[Necessity of quadratic roughness]
Let $x\in C^\alpha([0,T],\mathbb{R}^d)\cap Q_{\mathbb{T}}([0,T],\mathbb{R}^d)$. 
Let
$\pi=(\pi^n)_{n\geq 1}$ be a balanced partition sequence of $[0,T]$ such that $x\in Q_{\pi}([0,T],\mathbb{R}^d)$. Then:\\
$$\bigg( \forall t\in [0,T]),\  [x]_{\pi}(t) = [x]_{\mathbb{T}}(t)\ \bigg) \qquad \Rightarrow \qquad \forall  \beta \in (0,2\alpha),\; x\in R^\beta_\pi([0,T],\mathbb{R}^d).$$\label{lemma.stability}
\end{lemma}
\begin{proof}  
Take $\beta\in (0,2\alpha)$. Since $\pi$ is a balanced sequence of partitions and since for dyadic partition $\frac{|\mathbb{T}^{n}|}{|\mathbb{T}^{n+1}|}=2<\infty$, we can construct a subsequence $\left(\mathbb{T}^{l_n}\right)$ of $\left(\mathbb{T}^n \right)$ such that $ |\mathbb{T}^{l_n}| =O(|\pi^n|)$ as following.

\par Firstly, if $|\mathbb{T}^n|^\beta = O(|\pi^n|)$, then take $l_n=n$. Otherwise, for $\limsup_n  \frac{|\mathbb{T}^n|^\beta }{|\pi^n|}=\infty$,  we construct $l_n$ as follows,
\[l_n=\inf \{l\geq n: |\pi^n|\geq |\mathbb{T}^{l}|^{\beta}\}<\infty, \qquad \text{since } |\mathbb{T}^l|^\beta=\frac{1}{2^{l\beta}}\to 0.\label{eq.kkn} \]
Since $l_n\geq n$, we also have $l_n \to \infty$. So from the construction of $l_n$ we get the following inequality:
\[\forall n\geq 1, \qquad |\mathbb{T}^{l_n}|\leq |\pi^n|^{1/\beta} < |\mathbb{T}^{l_n-1}|.\]
Since the subsequence $(\mathbb{T}^{l_n})$ is also balanced, there exists constants $c_1$ and $c_2$ such that,
\[\forall n\geq 1, \qquad c_1N(\mathbb{T}^{l_n})\geq N(\pi^{n})^{1/\beta} > c_2N(\mathbb{T}^{l_n-1}).\] 
The points of the partition $\mathbb{T}^{l_n}$ are interspersed among those of $\pi^n$. Define for $k= 1, \cdots, N(\pi^n)$:
$$ p(n,k) = \inf\{ m\geq 1: \quad s^{n}_{m}\in (t^n_{k},t^n_{k+1}] \},$$ 
where, $\pi^n=(0=t^n_1<t^n_2<\cdots< t^n_{N(\pi^n)}=T)$ and $\mathbb{T}^{l_n} = (0=s^n_1<s^n_2<\cdots<s^n_{N(\mathbb{T}^{l_n})}=T)$. Then we have 
\begin{equation} s^{n}_{p(n,k)-1} \leq t^n_k < s^{n}_{p(n,k)}  <  \cdots < s^{n}_{p(n,k+1)-1} \leq t^n_{k+1} <  s^{n}_{p(n,k+1)}, \label{eq.order}\end{equation}
where, $p(n,N(\pi^n))-1 = N(\mathbb{T}^{l_n})$. 
From the construction of $l_n$ and the fact that $\limsup_n\frac{|\mathbb{T}^n|}{|\mathbb{T}^{n+1}|} = 2<\infty$, we can conclude that $|\pi^n|\asymp |\mathbb{T}^{l_n}|^\beta$. To prove this, assume for contradiction $|\pi^n|$ and $|\mathbb{T}^{l_n}|^\beta$ are not asymptotically comparable. Then  $\sup_{i=1,\cdots,N(\pi^n)}| p(n,i+1)- p(n,i)| \rightarrow \infty$ as $n\to\infty$. 
Then, from the above definition of $l_n$ and using the balanced property of $\pi^n $ and $\mathbb{T}^{l_n}$ we have $\limsup_n \frac{|\pi^n|}{|\mathbb{T}^{l_n}|^\beta}\rightarrow \infty $. Since $\limsup_n \frac{|\pi^n|}{|\mathbb{T}^{l_n-1}|^\beta} $ is bounded by $1$ for all $n\geq 1$, from the construction of $l_n$, we thus have $$\infty = \limsup_n\left(\frac{|\pi^n|}{|\mathbb{T}^{l_n}|^\beta}\right)-\limsup_n\left(\frac{|\pi^n|}{|\mathbb{T}^{l_n-1}|^\beta} \right)\leq \limsup_n \left(\frac{|\pi^n|}{|\mathbb{T}^{l_n}|^\beta}-\frac{|\pi^n|}{|\mathbb{T}^{l_n-1}|^\beta} \right)$$ 
$$=\limsup_n \frac{|\pi^n|}{|\mathbb{T}^{l_n-1}|^\beta}\left(\frac{|\mathbb{T}^{l_n-1}|^\beta}{|\mathbb{T}^{l_n}|^\beta} -1 \right)< \infty $$ 
which is a contradiction, and the last inequality is form the fact that $\frac{|\mathbb{T}^{l_n-1}|^\beta}{|\mathbb{T}^{l_n}|^\beta} =2^\beta<\infty$ and $ \limsup_n \frac{|\pi^n|}{|\mathbb{T}^{l_n-1}|^\beta} \leq 1$. Hence  the sequence $\sup_{i=1,\cdots,N(\pi^n)} \ p(n,i+1)- p(n,i)$ is bounded as  $n\to\infty$:
$$\exists M>0, \text{ such that } \forall i,n\geq 1,\qquad p(n,i+1)- p(n,i)\leq M< \infty.$$
Therefore $(\mathbb{T}^{l_n})^\beta_{n\geq 1}$ and $(\pi^{n})_{n\geq 1}$ are (asymptotically) comparable i.e. the sequences
\begin{equation}
    \frac{N(\mathbb{T}^{l_n})^{\beta}}{N(\pi^n)}  \qquad {\rm and}\qquad  \frac{|\mathbb{T}^{l_n}|^\beta}{|\pi^n|} 
\end{equation}
are uniformly bounded. Now since, $ \forall t\in [0,T],\  [x]_{\pi}(t) = [x]_{\mathbb{T}}(t)$, we have:
\[   \Bigg|[x]_{\pi^n}(t) - [x]_{\mathbb{T}^{l_n}}(t)\Bigg| \to 0. \]
For convenience denote $(\mathbb{T}^{l_n})= (d^n)$. We will only give the proof for $t=T$, for $t<T$ we will get one additional boundary term which goes to zero.
\par Decomposing $ \Delta^n_k=x(t^n_{k+1})-x(t^n_{k})$ along the partition points of $d^{n},$ we obtain,
$$
\underbrace{x(t^n_{k+1})-x(t^n_{k})}_{\Delta^n_k}= \underbrace{(x(t^n_{k+1})-x(s^{n}_{p(n,k+1)})}_{D_k} -  
\underbrace{( x(t^n_{k})-x(s^{n}_{p(n,k)}))}_{B_k} 
+ \underbrace{ \sum_{i=p(n,k)}^{p(n,k+1)-1} (x(s^{n}_{i+1})-x(s^{n}_{i}) )}_{C_k}.$$
 Grouping together the terms in $[x]_{d^{n}}$ according to Equation \eqref{eq.order} yields
$$[x]_{\pi^{n}}(T)-[x]_{d^{n}}(T)=\sum_{k=1}^{N(\pi^n)-1} \left[{\Delta^n_k}^t{\Delta^n_k}- \sum_{i=p(n,k)}^{p(n,k+1)-1}  \big(x(s^{n}_{i+1})-x(s^{n}_{i})\big)^t \big(x(s^{n}_{i+1})-x(s^{n}_{i})\big) \right]$$
$$=\sum_{k=1}^{N(\pi^n)-1} \left[{\Delta^n_k}^t{\Delta^n_k}- C_k^tC_k \right] + \sum_{k=1}^{N(\pi^n)-1} \left[C_k^tC_k -\sum_{i=p(n,k)}^{p(n,k+1)-1}  \big(x(s^{n}_{i+1})-x(s^{n}_{i})\big)^t \big(x(s^{n}_{i+1})-x(s^{n}_{i})\big)\right].$$
Now, the second term of the previous equation can also be represented as follows.
$$\sum_{k=1}^{N(\pi^n)-1}\left[C_k^tC_k -\sum_{i=p(n,k)}^{p(n,k+1)-1}  \big(x(s^{n}_{i+1})-x(s^{n}_{i})\big)^t \big(x(s^{n}_{i+1})-x(s^{n}_{i})\big)  \right]$$
$$=\sum_{k=1}^{N(\pi^n)-1} \bigg[ \left(x(s^{n}_{p(n,k+1)}) - x(s^{n}_{p(n,k)})\right)^t \left(x(s^{n}_{p(n,k+1)}) - x(s^{n}_{p(n,k)})\right) $$
$$ - \sum_{i=p(n,k)}^{p(n,k+1)-1}  \left(x(s^{n}_{i+1})-x(s^{n}_{i})\big)^t \big(x(s^{n}_{i+1})-x(s^{n}_{i})\right) \bigg] $$
$$= \sum_{j=1}^{N(\pi^n)-1}\sum_{s^n_i\neq s^n_{i'}\in (t^n_j,t^n_{j+1}]} \left(x(s^n_{i+1})-x(s^n_{i})\right)^t\left(x(s^n_{i'+1})-x(s^n_{i'})\right).$$
This is precisely the roughness term of $x$ along partition $\pi$. Hence to show that \[\sum_{j=1}^{N(\pi^n)-1}\sum_{t^n_i\neq t^n_{i'}\in (s^n_j,s^n_{j+1}]} \left(x(t^n_{i}\wedge t)-x(t^n_{i-1}\wedge t)\right)^t\left(x(t^n_{i'}\wedge t)-x(t^n_{i'-1}\wedge t)\right)\xrightarrow[]{n\to\infty} 0,\] 
\[\hspace*{-2cm}\text {we only have to show that: } \;\left|[x]_{\pi^{n}}(T)-[x]_{d^{n}}(T) - \left(\sum_{k=1}^{N(\pi^n)-1} {\Delta^n_k}^t{\Delta^n_k}- C_k^tC_k \right) \right|\to 0.\]
Now, From the assumption $\left|[x]_{\pi^{n}}(T)-[x]_{d^{n}}(T)\right| \to 0$. So we only need to show that\\  $\left| \sum_{k=1}^{N(\pi^n)-1} {\Delta^n_k}^t{\Delta^n_k}- C_k^tC_k  \right|\to 0.$ This is a consequence from the fact that quadratic variation along $\pi$ exists and $x\in C^{\alpha}([0,T],\mathbb{R}^d)$ with $\frac{\alpha}{\beta}>\frac{1}{2}$. A similar line of proof is adopted in more  detail in Theorem \ref{main.theorem}. 
\end{proof}
We will now show that quadratic roughness is also a {\it sufficient} condition for the uniqueness of quadratic variation along  balanced partition sequences.

Our main result is that quadratic roughness along such a sequence of partitions implies uniqueness of pathwise quadratic variation:
\begin{theorem}\label{main.theorem}
Let $\pi$ be a balanced sequence of partitions of $[0,T]$ and $x\in C^{\alpha}([0,T],\mathbb{R}^d)\cap R^\beta_\pi([0,T],\mathbb{R}^d)$ for some $0<\beta< 2\alpha$.
Then 
\[x\in Q_{\pi}([0,T],\mathbb{R}^d),\qquad \qquad \text{ and }\qquad \qquad \forall  t\in [0,T],\quad [x]_\pi(t)= [x]_{\mathbb{T}}(t).\]
\end{theorem}
\begin{proof}
Let, $\pi^n= \left(0=t^n_1<t^n_2<\cdots< t^n_{N(\pi^n)}=T\right)$.
Since $x\in R^\beta_\pi([0,T],\mathbb{R}^d)$, from Definition \ref{def.rough} we know there exists a sub/super-sequence $d= (d^n)_{n\geq 1}$ of the dyadic partition $\mathbb{T}= (\mathbb{T}^n)$ with  $d^n = \left(0=s^n_1<s^n_2<\cdots< s^n_{N(d^n)}=T\right)$ such that $ |d^n|^\beta = O(|\pi^n|)$ and for all $t\in [0,T]$ :
\begin{equation}
    \sum_{j=1}^{N(\pi^n)-1}\sum_{s^n_i\neq s^n_{i'}\in (t^n_j,t^n_{j+1}]} \left(x(s^n_{i+1}\wedge t)-x(s^n_{i}\wedge t)\right)^t\left(x(s^n_{i'+1}\wedge t)-x(s^n_{i'}\wedge t)\right)\xrightarrow[]{n\to\infty} 0.\label{rough.eq.maintheorem}
\end{equation}
So there exists $C< \infty$  and $ N_0\in \mathbb{N}$ such that:
\[ \forall n\geq N_0,  \qquad |d^n|^\beta\leq C |\pi^n|\]
We will assume $n\geq N_0$ for the rest of the proof.
Since both $d=(d^n)$ and $\pi=(\pi^n)$ are balanced sequence of partitions there exists $C_1<\infty$ such that:
\[ \forall n\geq 1,  \qquad  N(\pi^n)\leq C_1N(d^n)^\beta. \] 
 $d= (d^n)_{n\geq 1}$ is also a balanced sequence of partitions of $[0,T]$ and the points of  $\pi^n$ are interspersed along those of $d^n$. Define for $k= 1,2,\cdots, N(\pi^n)$:
$$ p(n,k) = \inf\{ m\geq 1: \quad s^{n}_{m}\in (t^n_{k},t^n_{k+1}] \}.$$ 
Then we get the following inequality regarding the partition points of $\pi^n$ and $d^n$. 
\begin{equation} s^{n}_{p(n,k)-1} \leq t^n_k < s^{n}_{p(n,k)}  <  \cdots < s^{n}_{p(n,k+1)-1} \leq t^n_{k+1} <  s^{n}_{p(n,k+1)}, \label{eq.ordering1}\end{equation}
where $p(n,N(\pi^n))-1 = N(d^{n})$ and $p(n,0) = 1 $. 
For all $t\in [0,T]$ we will show that by grouping the points of $d^{n}$ according to the intervals defined by $\pi^{n}$ and use the roughness property of $x$ along $\pi$:
$$[x]_{\pi^{n}}(t) =\sum_{k=1}^{N(\pi^n)-1} \big(x(t^n_{k+1}\wedge t)-x(t^n_{k}\wedge t)\big)^t \big(x(t^n_{k+1}\wedge t)-x(t^n_{k}\wedge t)\big)\quad{\rm and}$$
$$ [x]_{d^{n}}(t) =\sum_{k=1}^{N(d^{n})-1} \big(x(s^{n}_{k+1}\wedge t)-x(s^{n}_{k}\wedge t)\big)^t \big(x(s^{n}_{k+1}\wedge t)-x(s^{n}_{k}\wedge t)\big)$$
have the same limits. We define an auxiliary partition 
\[\sigma^n =  (0=s^n_{p(n,1)}<s^n_{p(n,2)}<\cdots <s^n_{p(n,N(\pi^n))-1} =T) \qquad \text{ and, }  \]
\[[x]_{\sigma^n}(t) =  \sum_{k=1}^{N(\pi^n)-1} \big(x(s^{n}_{p(n,k+1)}\wedge t)-x(s^{n}_{p(n,k)}\wedge t)\big)^t \big(x(s^{n}_{p(n,k+1)}\wedge t)-x(s^{n}_{p(n,k)}\wedge t)\big) \] 
and we will show that $[x]_{\pi^{n}}(t)$ and $[x]_{d^{n}}(t)$ have the same limit as $[x]_{\sigma^{n}}(t)$.
We shall give the proof for $t=T$; for $t<T$ we have an additional boundary term that goes to zero.

Decomposing $ \Delta^n_k=x(t^n_{k+1})-x(t^n_{k})$ along the partition points of $d^{n},$ we obtain,
$$
\underbrace{x(t^n_{k+1})-x(t^n_{k})}_{\Delta^n_k}= \underbrace{(x(t^n_{k+1})-x(s^{n}_{p(n,k+1)})}_{D_k} -  
\underbrace{( x(t^n_{k})-x(s^{n}_{p(n,k)}))}_{B_k} 
+ \underbrace{ \sum_{i=p(n,k)}^{p(n,k+1)-1} (x(s^{n}_{i+1})-x(s^{n}_{i}) )}_{C_k}.$$
 Grouping together the terms in $[x]_{d^{n}}$ according to Equation \eqref{eq.ordering1} yields
$$[x]_{\pi^{n}}(T)-[x]_{d^{n}}(T)=\sum_{k=1}^{N(\pi^n)-1} \left[{\Delta^n_k}^t{\Delta^n_k}- \sum_{i=p(n,k)}^{p(n,k+1)-1}  \big(x(s^{n}_{i+1})-x(s^{n}_{i})\big)^t \big(x(s^{n}_{i+1})-x(s^{n}_{i})\big) \right]$$
$$=\sum_{k=1}^{N(\pi^n)-1} \left[{\Delta^n_k}^t{\Delta^n_k}- C_k^tC_k \right] + \sum_{k=1}^{N(\pi^n)-1} \left[C_k^t C_k -\sum_{i=p(n,k)}^{p(n,k+1)-1}  \big(x(s^{n}_{i+1})-x(s^{n}_{i})\big)^t \big(x(s^{n}_{i+1})-x(s^{n}_{i})\big)  \right].$$
Now, the second term of the previous equation:
$$\sum_{k=1}^{N(\pi^n)-1}\left[C_k^tC_k -\sum_{i=p(n,k)}^{p(n,k+1)-1}  \big(x(s^{n}_{i+1})-x(s^{n}_{i})\big)^t \big(x(s^{n}_{i+1})-x(s^{n}_{i})\big)  \right]$$
$$=\sum_{k=1}^{N(\pi^n)-1} \bigg[ \left(x(s^{n}_{p(n,k+1)}) - x(s^{n}_{p(n,k)})\right)^t \left(x(s^{n}_{p(n,k+1)}) - x(s^{n}_{p(n,k)})\right) $$
$$ - \sum_{i=p(n,k)}^{p(n,k+1)-1}  \left(x(s^{n}_{i+1})-x(s^{n}_{i})\right)^t \left(x(s^{n}_{i+1})-x(s^{n}_{i})\right) \bigg] $$
$$=\sum_{j=1}^{N(\pi^n)-1}\sum_{s^n_i\neq s^n_{i'}\in (t^n_j,t^n_{j+1}]} \left(x(s^n_{i+1})-x(s^n_{i})\right)^t\left(x(s^n_{i'+1})-x(s^n_{i'})\right)\xrightarrow[]{n\to\infty} 0.$$
The last limit is precisely  \eqref{rough.eq.maintheorem} which arises from the roughness of $x$ along $\pi$. Additionally note that 
\[ \forall n\geq 1, \;[x]_{\sigma^n} = \sum_{k=0}^{N(\pi^n)-1}C_k^tC_k = [x]_{d^n}\quad \text{ and, }\quad\]\[ [x]_{\sigma} = \lim_{n\rightarrow\infty} [x]_{\sigma^n} = \lim_{n\rightarrow\infty} \sum_{k=0}^{N(\pi^n)-1}
C_k^tC_k = [x]_{d} = [x]_{\mathbb{T}} <\infty .\]
So $x\in Q_{\sigma}([0,T],\mathbb{R}^d)$.
Now to show that $|[x]_{\pi^{n}}(T)-[x]_{d^{n}}(T)| \rightarrow 0$ we only need to show that $\left|\sum_{k=1}^{N(\pi^n)-1} \left[{\Delta^n_k}^t{\Delta^n_k}- C_k^tC_k \right]\right|\rightarrow 0.$ Since,
$$\sum_{k=1}^{N(\pi^n)-1} \Delta_k^{nt} \Delta_k^n = \sum_{k=1}^{N(\pi^n)} ({C_k}+D_k-B_k)^t ({C_k}+D_k-B_k) $$
$$=\sum_{k=1}^{N(\pi^n)-1} {C_k}^t{C_k} +\sum_{k=1}^{N(\pi^n)-1} (D_k-B_k)^t(D_k-B_k) -2\sum_{k=1}^{N(\pi^n)-1} {C_k}^tB_k  +2\sum_{k=1}^{N(\pi^n)-1} {C_k}^tD_k.$$
we finally obtain,
$$\left|\sum_{k=1}^{N(\pi^n)-1}\left[{\Delta^n_k}^t{\Delta^n_k}-C_k^tC_k\right] \right|\leq \left|\sum_{k=1}^{N(\pi^n)-1} (D_k-B_k)^t(D_k-B_k)\right| +\left|2\sum_{k=1}^{N(\pi^n)-1} {C_k}^tB_k\right| +\left|2 \sum_{k=1}^{N(\pi^n)-1} {C_k}^tD_k\right|. $$
Now we will show that as $n\rightarrow \infty, \quad \left|\sum_{k=1}^{N(\pi^n)-1}({\Delta^n_k}^t{\Delta^n_k}-C_k^tC_k) \right| \rightarrow 0. $
Since $x \in C^\alpha([0,T],\mathbb{R}^d) $  we have :
$$ \forall t\in [0,T-h], \quad\forall h>0,\qquad \|x(t+h)-x(t)\| \leq \|x\|_\alpha h^{\alpha}.$$
Now, 
$$\left|\sum_{k=1}^{N(\pi^n)-1}D_k^tD_k\right|\leq \sum_{k=1}^{N(\pi^n)-1}||D_k||^2 \leq \sum_{k=1}^{N(\pi^n)-1} \|x\|_\alpha^2|d^{n}|^{2\alpha}$$ 
$$ \leq \|x\|_\alpha^2 N(\pi^n) |d^{n}|^{2\alpha}\leq cC^{\frac{1}{\beta}} N(\pi^n) |\pi^{n}|^{2\alpha/\beta}\mathop{\to}^{n\to \infty} 0$$
since $\frac{2\alpha}{\beta}>1$. Similarly we have $\sum_{k=1}^{N(\pi^n)-1} B_k^t B_k \rightarrow 0$. Therefore, 
$$\sum_{k=1}^{N(\pi^n)-1}|(D_k-B_k)^t (D_k-B_k)| \leq 2\left|\sum_{k=1}^{N(\pi^n)-1} D_k^tD_k \right| + 2 \left|\sum_{k=1}^{N(\pi^n)-1} B_k^tB_k \right| \rightarrow 0.$$
Using  H\"older's inequality,
$$   \left|\sum_{k=1}^{N(\pi^n)-1} D_k^t{C_k} \right| \leq   \left(\sum_{k=1}^{N(\pi^n)-1} \|D_k\|^2\right)^{\frac{1}{2}} \left( \sum_{k=1}^{N(\pi^n)-1} \|{C_k}\|^2 \right)^{\frac{1}{2}}.$$
Since the quadratic variation of $x$ along the sequence of partitions $\sigma$ exists and finite;  the sequence $ \sum_{k=1}^{N(\pi^n)-1} \|{C_k}\|^2$ is bounded. Combining this with the estimate above we obtain $$\left|\sum_{k=1}^{N(\pi^n)-1} D_k^t{C_k}\right| \rightarrow 0.$$
Similarly, we have, $\left|\sum_{k=1}^{N(\pi^n)-1} B_k^t{C_k} \right|\rightarrow 0$ as $n \to \infty$.
Therefore, $\sum_{k=1}^{N(\pi^n)-1}\left[{\Delta^n_k}^t{\Delta^n_k}-C_k^tC_k\right] \rightarrow 0$. Hence,
$$\Big|[x]_{\pi^{n}}(T)-[x]_{d^{n}}(T)\Big|$$
$$=\left|\sum_{k=1}^{N(\pi^n)-1} \left[{\Delta^n_k}^t{\Delta^n_k}- C_k^tC_k \right] + \sum_{k=1}^{N(\pi^n)-1} \left[C_k^tC_k -\sum_{i=p(n,k)}^{p(n,k+1)-1}  \big(x(s^{n}_{i})-x(s^{n}_{i-1})\big)^t \big(x(s^{n}_{i})-x(s^{n}_{i-1})\big) \right]\right|$$
$$\leq \left|\sum_{k=1}^{N(\pi^n)-1} \left[{\Delta^n_k}^t{\Delta^n_k}- C_k^tC_k \right] \right|+\left| \sum_{k=1}^{N(\pi^n)-1} \left[C_k^tC_k -\sum_{i=p(n,k)}^{p(n,k+1)-1}  \big(x(s^{n}_{i})-x(s^{n}_{i-1})\big)^t \big(x(s^{n}_{i})-x(s^{n}_{i-1})\big) \right]\right|\xrightarrow[]{n\to \infty}0.$$
Since $d=(d^n)_{n\geq 1}$ is a sub/super-sequence of the dyadic partition, for all $t\in [0,T]$, $[x]_{d}= [x]_{\mathbb{T}}$. So, $x\in Q_{\pi}([0,T],\mathbb{R}^d)$ and also, $$\forall t \in [0,T], \qquad [x]_{\pi}(t)= [x]_{\mathbb{T}}(t).$$
This concludes the proof.
\end{proof} 

For H\"older continuous  paths, the quadratic roughness property along a balanced partition sequence implies existence of quadratic variation along the same sequence:
\begin{corollary}
Let  $\pi= (\pi^n)_{n\geq 1}$ be a balanced partition sequence of $[0,T]$ with $|\pi^n| \to 0 $. Then 
$$\forall \beta\in(0,2\alpha), \qquad R^\beta_\pi([0,T],\mathbb{R}^d)\cap C^\alpha([0,T],\mathbb{R}^d) \subset  Q_{\pi}([0,T],\mathbb{R}^d)\cap Q_{\mathbb{T}}([0,T],\mathbb{R}^d).$$
\end{corollary}
\begin{proof}
Let $\pi\in \mathbb{B}([0,T])$. If $x\in R^\beta_\pi([0,T],\mathbb{R}^d)\cap C^\alpha([0,T],\mathbb{R}^d)$ then from  Definition \ref{def.rough}, $x\in Q_{\mathbb{T}}([0,T],\mathbb{R}^d)$. Since $\beta\in (0,2\alpha)$, from Theorem \ref{main.theorem} the quadratic variation of $x $ along $\pi $ exists and is equal its quadratic variation along the dyadic partition. So $x\in Q_{\pi}([0,T],\mathbb{R}^d)$.
\end{proof}
In general without the H\"older continuity assumption on  $x\in C^0([0, T],\mathbb{R}^d)$, roughness along a partition sequence $\pi$ does not imply the existence of quadratic variation along  $\pi$. 
The following lemma is a simple application of Theorem \ref{main.theorem}:
\begin{lemma}
Let $\pi=(\pi^n)_{n\geq 1}$ and $\sigma=(\sigma^n)_{n\geq 1}$ be balanced sequences of partitions of $[0,T]$. If $x\in C^\alpha([0,T],\mathbb{R}^d)\cap R^\beta_\pi([0,T],\mathbb{R}^d)\cap R^\gamma_\sigma([0,T],\mathbb{R}^d)$ for some $\beta,\gamma \in (0,2\alpha)$ then: 
\[ x\in Q_{\pi}([0,T],\mathbb{R}^d)\cap Q_{\sigma}([0,T],\mathbb{R}^d)\qquad \text{ and, } \qquad\forall t\in [0,T] \qquad [x]_{\pi}(t) = [x]_{\sigma}(t).\]
\end{lemma}
 
\begin{corollary}
Let $\pi=(\pi^n)_{n\geq 1}$  be a balanced sequences of partitions of $[0,T]$. If $x\in R^\beta_\pi([0,T],\mathbb{R}) \cap C^\alpha([0,T],\mathbb{R})$ for some $0< \beta< 2\alpha\wedge 1$ then:
\[\forall f\in C^2([0,T],\mathbb{R}), \forall t\in[0,T]: \quad [f\circ x]_\pi(t) = [f\circ x]_\mathbb{T}(t)\]
\end{corollary}
\begin{proof}
Since $\pi$ and $x$ satisfies the conditions of Theorem \ref{main.theorem}, we can conclude $\forall t\in [0,T]: \; [x]_\pi(t) = [x]_\mathbb{T}(t)$. So $[f\circ x]_\pi(t)$ can be expresses as:
\[[f\circ x]_\pi(t)=\lim_{n\to \infty} \sum_{\pi^n\cap [0,t]} \left(f\circ x(t^n_{i+1})- f\circ x(t^n_i)\right)^2 = \int_0^t f'\circ x(u) d[x]_\pi(u) \]
\[= \int_0^t f'\circ x(u) d[x]_\mathbb{T}(u)   =\lim_{n\to \infty} \sum_{\mathbb{T}^n\cap [0,t]} \left(f\circ x(s^n_{j+1})- f\circ x(s^n_j)\right)^2 = [f\circ x]_{\mathbb{T}}(t).\]
\end{proof}

\subsection{Invariant definition of quadratic variation} 
Let $\mathbb{T} =(\mathbb{T}^n)_{n\geq 1}$ be the dyadic sequence of partitions of $[0,T]$. Define,
\begin{equation}
   {\cal Q}([0,T],\mathbb{R}^d ) =  C^{\frac{1}{2} -}([0,T],\mathbb{R}^d )  \cap Q_{\mathbb{T}}([0,T],\mathbb{R}^d )  .\label{eq.Q}
\end{equation}

\begin{lemma}
The class ${\cal Q}([0,T],\mathbb{R}^d )$ is non-empty and contains all `typical' Brownian paths.
\end{lemma}

\begin{proof}
 Let $W$ be a Wiener process on a probability space $(\Omega, {\cal F},\mathbb{P})$, which we take to be the canonical Wiener space without loss of generality. For dyadic partition $\mathbb{T}$, since we have $|\mathbb{T}^n|= \frac{1}{2^n}$, so $|\mathbb{T}^n|\log (n)\rightarrow 0$. So from Dudley \cite{dudley2011} we can conclude:
 \[\mathbb{P}\left[W\in Q_{\mathbb{T}}\left([0,T],\mathbb{R}^d \right) \right]=1. \]
Brownian paths are almost-surely $\alpha$-H\"older   for  $\alpha <\frac{1}{2}$, so 
$$\mathbb{P}\left( W\in Q_{\mathbb{T}}([0,T],\mathbb{R}^d ) \cap C^{\frac{1}{2}-}([0,T],\mathbb{R}^d) \right)=1, $$
hence the result follows.
\end{proof}
Based on the results above we can now give an `intrinsic' definition of pathwise quadratic variation  for paths in ${\cal Q}([0,T],\mathbb{R}^d )$ which does not rely on a particular partition sequence:
\begin{proposition}[Quadratic variation map] \label{prop.qv} There exists a unique map:
\begin{eqnarray} [\ .\ ]\quad :\ {\cal Q}\left([0,T],\mathbb{R}^d \right) & \to  & {C}^0\left([0,T],S^+_d \right)\nonumber
\\x\ & \mapsto & [x]\nonumber
\end{eqnarray}
such that
$$\forall \pi \in \mathbb{B}([0,T]),\; \forall \beta\in(0,1),\quad \forall x\in {  R}^\beta_\pi([0,T],\mathbb{R}^d )\cap {\cal Q}([0,T],\mathbb{R}^d ),\quad \forall t\in [0,T],\quad [x]_\pi(t)=[x](t).$$
We call $[x]$ the quadratic variation of $x$.
\end{proposition}
\begin{proof}
Let $\pi\in \mathbb{B}([0,T])$. Then for $\beta\in (0,1)$ and for any $x\in {  R}^\beta_{\pi}([0,T],\mathbb{R}^d )\cap {\cal Q}([0,T],\mathbb{R}^d )$ take $\alpha= \frac{\beta+1}{4}<\frac{1}{2}$ as $\beta<1$. So we have 
\[x\in  {  R}^\beta_{\pi}([0,T],\mathbb{R}^d )\cap  {  C}^\alpha([0,T],\mathbb{R}^d ) \cap {Q}_{\mathbb{T}}([0,T],\mathbb{R}^d ). \]
Then  Theorem \ref{main.theorem} implies 
$$  x\in Q_{\pi}([0,T],\mathbb{R}^d)\qquad \text{ and, }\qquad \forall t\in [0,T],\; [x]_{\pi}(t) = [x]_{\mathbb{T}}(t). $$
By the same argument the quadratic variation does not depend  on the choice of  $\pi\in \mathbb{B}([0, T])$, so the result follows.
\end{proof}
\begin{remark} If $X$ is a continuous $\mathbb{P}$-semimartingale then its image $[X]$ under the map defined in Proposition \ref{prop.qv} coincides almost-surely with the probabilistic definition of quadratic variation as a limit in probability \cite{karandikar1983,protter}. Building on \cite{karandikar1983}, Karandikar and Rao \cite{karandikar2014} construct a (different) quadratic variation map which shares this property. In contrast to  \cite{karandikar2014}, our construction does not use any probabilistic tools, does not rely on specific path-dependent partitions and identifies explicitly the domain of definition of the map (rather than implicitly in terms of the support of a probability measure).
\end{remark}

\section{Pathwise It\^o calculus}
\subsection{Pathwise integration and the \follmer-It\^o formula}\label{sec:ito}

Using Theorem \ref{main.theorem} and Proposition \ref{prop.qv}, we can give a formulation of  \follmer's pathwise It\^o calculus which is invariant with respect to the choice of the  sequence of partitions $\pi$.
\begin{theorem}[Invariance of the F\"ollmer integral]
 There exists a unique map
\begin{eqnarray} I\quad :\ C^2(\mathbb{R}^d)\times{\cal Q}([0,T],\mathbb{R}^d ) & \to  & {\cal Q}([0,T],\mathbb{R} )\nonumber
\\ (f, x) \ & \to & I(f,x)=\int_0^. (\nabla f\circ x).dx,\nonumber
\end{eqnarray}
such that:
$\forall \pi \in \mathbb{B}([0,T]),\quad \forall \beta\in (0,1), \quad \forall x\in {  R}^\beta_\pi([0,T],\mathbb{R}^d )\cap {\cal Q}([0,T],\mathbb{R}^d ),\quad \forall t\in [0,T],$
$$I(f,x)(t)=\int_0^t (\nabla f\circ x).d^\pi x=\mathop{\lim}_{n\to\infty} \sum_{\pi^n} \nabla f(x(t^n_i)).(x(t^n_{i+1}\wedge t)-x(t^n_i\wedge t) ).$$
We denote $I(f,x)=\int_0^. (\nabla f\circ x) dx$.
Furthermore, we have the following change of variable formula. 
\ba
\forall f\in C^2(\mathbb{R}^d),\quad  \forall \pi\in \mathbb{B} ([0,T]),\quad \forall\beta\in (0,1) ,\quad \forall x\in {  R}^{\beta}_\pi([0,T],\mathbb{R}^d )\cap {\cal Q}([0,T],\mathbb{R}^d ),\; \forall t\in [0,T]\nonumber\\
\label{int.ind.parti} f(x(t))-f(x(0))= \int_0^t (\nabla f\circ x).dx +\frac{1}{2}\int_0^t <\nabla^2f(x),d[x]>\qquad\qquad\qquad\\
{\rm and}\qquad \left[\int_0^. (\nabla f\circ x)\; dx)\right]_{\pi}(t)=\int_0^t <(\nabla f\circ x)^t(\nabla f\circ x),d[x]>.\qquad\qquad\qquad\ea
\label{thm.ito} 
\end{theorem}
\begin{proof}
For any $\beta\in (0,1)$, take $\alpha  = \frac{\beta+1}{4}<\frac{1}{2}$. Fix any balanced partition sequence $\pi$, if $x\in {  R}^{\beta}_\pi([0,T],\mathbb{R}^d )\cap {\cal Q}([0,T],\mathbb{R}^d )$ then using Theorem \ref{main.theorem} we can conclude $x\in Q_\pi([0,T],\mathbb{R}^d)\cap Q_\mathbb{T}([0,T],\mathbb{R}^d)$. So for the balanced partition sequence $\pi\in \mathbb{B}([0,T])$ the pathwise It\^o formula \cite{follmer1981} implies
$$\int_0^t (\nabla f\circ x).d^{\pi}x= f(x(t))-f(x(0)) -\frac{1}{2}\int_0^t <\nabla^2f(x),d[x]_{\pi}>,  \qquad{\rm and}$$
$$\int_0^t (\nabla f\circ x).d^{\mathbb{T}}x = f(x(t))-f(x(0))-\frac{1}{2}\int_0^t <\nabla^2f(x),d[x]_{\mathbb{T}}>.$$
Since all  assumptions of Theorem \ref{main.theorem} are satisfied for the path $x$ along the sequence of partitions $\pi$, we conclude 
 $[x]_{\mathbb{T}} = [x]_{\pi} $. This argument is true for all $\beta\in(0,1)$ and for all $\pi\in \mathbb{B}([0,T])$. So:
 $$\forall \beta\in(0,1),\; \forall \pi\in \mathbb{B}([0,T]) \;\forall t\in [0,T]\;:\quad \int_0^t (\nabla f\circ x).d^{\pi}x= \int_0^t (\nabla f\circ x).d^{\mathbb{T}}x $$
 i.e. the pathwise integral $\int_0^t  (\nabla f\circ x).d^{\pi}x$ along a balanced sequence of partitions $\pi$ does not depend on choice of   $\pi$. 
To show $I(f,x)\in {\cal Q}([0,T],\mathbb{R}) $ we first note that
by \cite[Lemma 4.11]{ananova2018} we have $ I(f,x) \in C^{\frac{1}{2}-}([0,T],\mathbb{R}).$

For all $\beta\in (0,1)$ and for all $\pi \in \mathbb{B}([0,T])$ using the Theorem \ref{main.theorem} we can conclude $x\in {\cal Q}([0,T],\mathbb{R})\cap {\cal R}_{\pi}^{\beta}([0,T],\mathbb{R})$ implies $x\in {Q}_{\pi}([0,T],\mathbb{R})$. Now applying the pathwise isometry formula \cite[Theorem 2.1]{ananova2017}, to the integral $\int_0^. (\nabla f\circ x).d^{\pi}x$ we obtain that $\int_0^. (\nabla f\circ x).d^{\pi}x=\int_0^. (\nabla f\circ x).d^{\mathbb{T}}x\in Q_\pi([0,T],\mathbb{R})\cap Q_\mathbb{T}([0,T],\mathbb{R})$ and
$$ \left[\int_0^. \nabla f\circ x\; dx\right]_{\pi}(t)=\int_0^t <(\nabla f\circ x)^t(\nabla f\circ x),d[x]_{\pi}>. $$
From Theorem \ref{main.theorem} we have $[x]_{\mathbb{T}} = [x]_{\pi} $, so $ \left[\int_0^. \nabla f\circ x\; dx\right]_{\pi}(t) $
does not depend on choice of balanced partition $\pi$. As a consequence:
$$\left[\int_0^. \nabla f\circ x \; dx\right]_{\pi}(t) =  \left[\int_0^. \nabla f\circ x\; dx\right]_{\mathbb{T}}(t)=\int_0^t <(\nabla f\circ x)^t(\nabla f\circ x),d[x]>.$$
So finally $I(f,x)\in{ \cal Q}([0,T],\mathbb{R})$.
\end{proof}

\subsection{Local time}\label{sec:localtime}
Pathwise analogues of (semimartingale) local time have been considered in \cite{bertoin1987,perkowski2019,davis2014,kim2019,perkowski2015,wuermli1980} in the context of Tanaka-type  formulas for convex functions or functions with Sobolev regularity. One such construction of local time involves taking a limit of a sequence of discrete approximations  of occupation densities along a fixed sequence of time partitions \cite{perkowski2019,davis2018}.

Given a partition sequence  $\pi=(\pi^n)_{n\geq 1}$ and a path $x\in C^0([0,T],\mathbb{R})\cap Q_\pi([0,T],\mathbb{R})$, define the function
 $L^{\pi^n}_t : \mathbb{R}\to \mathbb{R}$ by
$$L^{\pi^n}_t (u) := 2 \sum_{t^n_j\in \pi^n\cap [0,t]} \mathbbm{1}_{[[x(t^n_{j}),x(t^n_{j+1}))]}(u)\;|x(t^n_{j+1}\wedge t) - u|. $$ 
where $[[u,v)]:= 
               [u,v)$ if $u\leq v$ and
               $[[u,v)]:=[v,u)$ if $u > v$.
            $L^{\pi^n}_t$ is bounded and zero outside $[\min x, \max x]$.
               
Following \cite{wuermli1980,bertoin1987,davis2014,perkowski2015}
we say that $x$ has  ($\mathcal{L}^2$)- local time on $[0,T]$ along $\pi$ if the sequence $(L^{\pi^n}_t,n \geq 1)$ converges weakly in $\mathcal{L}^2(\mathbb{R})$  to a limit  $ L^{\pi}_t$ for all $t\in [0,T]$:
 $$\forall t\in [0,T],\quad \forall h\in \mathcal{L}^2(\mathbb{R}), \quad \int L^{\pi^n}_t(u) h(u) du  \xrightarrow[]{n\to\infty} \int L^{\pi}_t(u) h(u) du.$$
The local time along $\pi$ satisfies the  occupation time formula \cite{wuermli1980,bertoin1987,perkowski2015}: for every Borel set $A \in \mathcal{B}(\mathbb{R})$,
$$   \int_{A} L^{\pi}_{t}(u)du= \frac{1}{2} \int_{0}^{t} \mathbbm{1}_{A}(x) d[x]_{\pi} $$
and the following extension of the pathwise Ito formula 
\eqref{int.ind.parti}
to functions in the Sobolev space $W^{2,2}(\mathbb{R})$ (see e.g.\cite[Thm 3.1]{davis2018}):  
\ba\forall f\in W^{2,2}(\mathbb{R}),\; \forall t\in [0,T],\quad f(x(t))-f(x(0))=\int_0^t (f'\circ x).d^\pi x+\frac{1}{2}\int_\mathbb{R}L^{\pi}_t(u) f''(u) du,\;\label{eq.tanaka}\ea
where the first integral is a limit of left Riemann sums along time partition $\pi$:
$$\int_0^t (f'\circ x).d^\pi x:= \lim_{n\to\infty} \sum_{\pi^n\cap [0,t]} f'(x(t^n_i)).(x(t^n_{i+1})-x(t^n_i)).$$
Unlike the intrinsic definition of occupation densities for real functions (see e.g. \cite{geman1976}), the above construction depends on the choice of the partition sequence $\pi$ and a natural question is therefore to clarify the dependence of this local time on the choice of the partition sequence.
Note that, differently from \cite{geman1976,geman1980}, $L^{\pi}_t$ is the density of a {\it weighted} occupation measure, weighted by quadratic variation $[x]_\pi$ so a necessary condition for the uniqueness of $L^{\pi}_t$ is the uniqueness of $[x]_\pi$.

We now show that the quadratic roughness property  implies an invariance property of
the { local time}  with respect to the sequence of partitions:
\begin{theorem}[Invariance of local time under quadratic roughness]\label{theorem.localtime}
Let $x\in C^{\alpha}([0,T],\mathbb{R})\cap$ $R^{\beta}_\pi([0,T],\mathbb{R})$ with $0< \beta\leq 2\alpha < 1$. Assume $x$ has local time $L^{\pi}_t$ on $[0,t]$ along $\pi\in \mathbb{B}(0,T])$.
Furthermore if $x$ also has local time $L^{\mathbb{T}}_{t}$ on $[0,t] $ along the dyadic partition sequence $\mathbb{T}$, then we can conclude 
\[\forall t\in [0,T], \qquad   L^{\pi}_{t}(u) = L^{\mathbb{T}}_{t}(u)\qquad du-a.e.\] This defines a unique  element $L_{t}\in \mathcal{L}^2(\mathbb{R})$ which we call the local time of $x$ on $[0,t]$.
\end{theorem}
This result shows that for  paths satisfying the quadratic roughness property, the ($\mathcal{L}^2$)- local time  is an intrinsic object associated with the path $x$, independent of the  (balanced) sequence of partitions used in the construction.
\begin{proof} 
From \cite[Satz 9]{wuermli1980} for any Borel set $A \in \mathcal{B}(\mathbb{R})$ we have the occupation density formula as: 
$$  \int_{A} L^{\pi}_{t}(u)du= \frac{1}{2} \int_{0}^{t} \mathbbm{1}_{A}(x) d[x]_{\pi}. $$
Since the local time along $\mathbb{T}$ exists, we also have: 
$$ \forall A \in \mathcal{B}(\mathbb{R}),\quad  \int_{A} L^{\mathbb{T}}_{t}(u)du= \frac{1}{2} \int_{0}^{t} \mathbbm{1}_{A}(x) d[x]_{\mathbb{T}}. $$
As $\pi$ is balanced, Theorem \ref{main.theorem} implies that $[x]_{\pi}= [x]_{\mathbb{T}}$. Hence, 
$$\forall A \in \mathcal{B}(\mathbb{R}),\quad \int_{A} L^{\pi}_{t}(u)du= \int_{A} L^{\mathbb{T}}_{t}(u)du,$$
which implies $L^\pi_t=L^\mathbb{T}_t$ almost everywhere.
\end{proof}

An important consequence of  this result is the uniqueness of limits of left Riemann sums for integrands in the Sobolev space $W^{1,2}(\mathbb{R})$ and a robust version of the pathwise Tanaka formula \cite{bertoin1987,davis2018}:
\begin{corollary}[Uniqueness of F\"ollmer integral on $W^{1,2}(\mathbb{R})$ and pathwise Tanaka formula]\ \\
Under the assumptions of Theorem \ref{theorem.localtime} we also have:
$$\forall h\in W^{1,2}(\mathbb{R}),\forall t\in [0,T],\qquad \int_0^t (h\circ x)d^\pi x= \int_0^t (h\circ x)d^\mathbb{T} x.$$
Designating this common value by $\int_0^t (h\circ x)dx$, we obtain
\ba\forall f\in W^{2,2}(\mathbb{R}),\forall t\in [0,T],\quad f(x(t))-f(x(0))=\int_0^t (f'\circ x).dx+\frac{1}{2}\int_\mathbb{R}L_t(u) f''(u) du,\label{eq.robustanaka}\ea
where the pathwise integral and the local time may be computed with respect to {\it any} balanced partition sequence along which $x$ has quadratic roughness and local time.
\end{corollary}

\begin{acks}[Acknowledgments]
The authors would like to thank Anna Ananova, Jan Ob{\l}{\'o}j and David Pr\"omel for constructive comments and helpful discussions.
\end{acks}


\appendix
\section*{Technical proofs}
\subsection{Proof of Proposition \ref{prop.well.balanced}}\label{appendix.2.2}

(i) For any sequence of partitions $\pi$ of $[0,T]$ and for any $n\geq 1$: 
\[N(\pi^n)\underline{\pi^n} 
\leq T \leq N(\pi^n)|\pi^n|.\]
For proof of ($\Rightarrow$): Using the balanced property, $$\liminf_{n\rightarrow \infty} N(\pi^n)\underline{\pi^n} = \liminf_{n\rightarrow \infty} N(\pi^n)|\pi^n| \frac{\underline{\pi^n}}{|\pi^n|} \geq \liminf_{n\rightarrow \infty} \frac{1}{c} N(\pi^n)|\pi^n|\geq \frac{T}{c} >0.$$
Similarly, $$\limsup_{n\rightarrow \infty} N(\pi^n)|\pi^n| = \limsup_{n\rightarrow \infty} N(\pi^n)\underline{\pi^n} \frac{|\pi^n|}{\underline{\pi^n}} \leq \limsup_{n\rightarrow \infty} c N(\pi^n)\underline{\pi^n} \leq cT < \infty.$$ 
\\For proof of ($\Leftarrow$): $\limsup_{n\rightarrow \infty} \frac{|\pi^n|}{\underline{\pi^n}} =\limsup_{n\rightarrow \infty} \frac{N(\pi^n)|\pi^n|}{N(\pi^n)\underline{\pi^n}} \leq \frac{\limsup_{n\rightarrow \infty}N(\pi^n)|\pi^n|}{\liminf_{n\rightarrow \infty}N(\pi^n)\underline{\pi^n}} < \infty. $\\
(ii) For any sequence of partitions $\pi$ with vanishing mesh and for any fixed $h>0$ there exists a $N_0$ such that for all $n\geq N_0, \; |\pi^n|<h $. So for all $n\geq N_0$ and for   $t\in[0,T-h],\; N(\pi^n,t,t+h)\geq 1$. Hence:
$$ \underline{\pi^n} \leq \frac{h}{N(\pi^n,t,t+h)}\leq |\pi^n|. $$
So
$$\limsup_{n\rightarrow \infty} \frac{\sup_{t\in[0,T-h]} N(\pi^n,t,t+h)}{\inf_{t\in[0,T-h]} N(\pi^n,t,t+h)} \leq \limsup_{n\rightarrow \infty} \frac{|\pi^n|}{h} \times\frac{h}{\underline{\pi^n}} <\infty. $$
\\(iii) For any balanced sequence of partitions $\pi$ of $[0,T]$ and for any $n\geq 1$: 
$$\frac{1}{c}N(\pi^n)|\pi^n|\leq N(\pi^n)\underline{\pi^n} \leq T \leq N(\pi^n)|\pi^n|\leq cN(\pi^n)\underline{\pi^n},$$
where, $c$ is a constant $< \infty$. So the equivalence follows.
\\(iv) Let $\pi = (\pi^n)_{n\geq 1}$ be any balanced sequence of partitions of $[0,T]$:
$$ \pi^n=(0=t^n_1<t^n_2<\cdots <t^n_{N(\pi^n)}=T). $$
Now, define the new partition $g(\pi)=(g(\pi^n))_{n\geq 1}$ on $[g(0),g(T)]$ as follows:
$$ g(\pi^n) = \left(g(0)=g(t^n_1)<g(t^n_2)<\cdots <g(t^n_{N(\pi^n)}=g(T))\right). $$
Now, from mean value theorem there exists $u^n_k,v^n_k\in[t^n_k,t^n_{k+1}]$ such that,

$$\left|\limsup_{n\rightarrow\infty} \frac{\left|g(\pi^n)\right|}{\underline{g(\pi^n)}}\right|=\left|\limsup_{n\rightarrow \infty} \frac{\sup_{\pi^n} (g(t^n_{k+1})-g(t^n_k))}{\inf_{\pi^n} (g(t^n_{k+1})-g(t^n_k))}\right| =\left| \limsup_{n\rightarrow \infty} \frac{\sup_{\pi^n} g'(u^n_k) (t^n_{k+1}-t^n_k)}{\inf_{\pi^n} g'(v^n_k) (t^n_{k+1}-t^n_k)}\right| $$
$$\leq \left| \limsup_{n\rightarrow \infty} \frac{\sup_{\pi^n} g'(u^n_k)}{\inf_{\pi^n} g'(v^n_k)}\right|\times \left| \limsup_{n\rightarrow \infty} \frac{\sup_{\pi^n} (t^n_{k+1}-t^n_k)}{\inf_{\pi^n} (t^n_{k+1}-t^n_k)}\right| \leq \frac{\sup g'}{\inf g'}\times c <\infty.$$

 \subsection{Proof of Proposition \ref{prop.rough.fn}}\label{appendix.A2}
The proof of $1-3$ are direct consequence of Definition \ref{def.rough}.
\\\textbf{Proof of $4$}:
We will proof the statement for $d=1$ first, and then give general argument for $d>1$.
Let $\pi^n = (0=t^n_0<t^n_2<\cdots< t^n_{N(\pi^n)}=T)$. Since $x\in R^\beta_{\pi}([0,T],\mathbb{R}^d)$, 
there exists a dyadic sub/super-sequence $(d^n)_{n\geq 1}$ satisfying (i) $|d^n|^{\beta} =O\left(|\pi^n|\right)$ and
(ii) For all $t\in [0,T]$ : \begin{equation}\label{eq.P}
\sum_{k=1}^{N(\pi^n)}\sum_{s^n_i\neq s^n_{i'}\in (t^n_k,t^n_{k+1}]} \left(x(s^n_{i+1}\wedge t)-x(s^n_{i}\wedge t)\right)^t\left(x(s^n_{i'+1}\wedge t)-x(s^n_{i'}\wedge t)\right)\xrightarrow[]{n\to\infty} 0,\end{equation}
where, $d^n= \left(0=t^n_1<t^n_2<\cdots<t^n_{N(d^n)}=T\right)$. Consider now
\[[x]_{\pi^{n}}(t) =\sum_{k=0}^{N(\pi^n)-1} \big(x(t^n_{k+1}\wedge t)-x(t^n_{k}\wedge t)\big)^t \big(x(t^n_{k+1}\wedge t)-x(t^n_{k}\wedge t)\big),\quad{\rm and}\]
\[[x]_{d^{n}}(t) =\sum_{k=0}^{N(d^{n})-1} \left(x(s^{n}_{k+1}\wedge t)-x(s^{n}_{k}\wedge t)\right)^t \left(x(s^{n}_{k+1}\wedge t)-x(s^{n}_{k}\wedge t)\right).\]
Define for $k= 1,2,\cdots, N(\pi^n)$,
 $p(n,k) = \inf\{ m\geq 1: \quad s^{n}_{m}\in (t^n_{k},t^n_{k+1}] \}.$\\
Then we have the following inequality between partition points of $\pi^n$ and $d^n$:
\begin{equation} s^{n}_{p(n,k)-1} \leq t^n_k < s^{n}_{p(n,k)} < \cdots < s^{n}_{p(n,k+1)-1} \leq t^n_{k+1} < s^{n}_{p(n,k+1)}, \label{eq.ordering2}\end{equation}
where $p(n,N(\pi^n))-1 = N(d^{n})$. We define an auxiliary partition $\sigma=(\sigma^n)$ and show that $[x]_{\pi^{n}}(t)$ and $[x]_{d^{n}}(t)$ has exact same limit with $[x]_{\sigma^n}(t)$:
\[\sigma^n = (0=s^n_{p(n,1)}<s^n_{p(n,1)}<\cdots <s^n_{p(n,N(\pi^n))-1} =T) \qquad \text{ and, } \]
\[[x]_{\sigma^n}(t) = \sum_{k=0}^{N(\pi^n)-1} \big(x(s^{n}_{p(n,k+1)}\wedge t)-x(s^{n}_{p(n,k)}\wedge t)\big)^t \big(x(s^{n}_{p(n,k+1)}\wedge t)-x(s^{n}_{p(n,k)}\wedge t)\big).\] 
The roughness assumption precisely tells us for all $t\in [0,T]$, we have $[x]_{\sigma}(t) = [x]_{d}(t) =[x]_{\mathbb{T}}(t)$. So $x\in Q_\sigma([0,T],\mathbb{R}^d)$.
Hence to show that $x+y\in R^\beta_\pi([0,T],\mathbb{R}^d)$ we have to show that $[x+y]_{\sigma}(t) = [x+y]_{d}(t) =[x+y]_{\mathbb{T}}(t)$.
Now, 
\[[x+y]_{\sigma^n}(t) = \sum_{k=0}^{N(\pi^n)-1} \big((x+y)(s^{n}_{p(n,k+1)}\wedge t)-(x+y)(s^{n}_{p(n,k)}\wedge t)\big)^t \big((x+y)(s^{n}_{p(n,k+1)}\wedge t)-(x+y)(s^{n}_{p(n,k)}\wedge t)\big)\]
\[ = \sum_{k=0}^{N(\pi^n)-1} \big(x(s^{n}_{p(n,k+1)}\wedge t)-x(s^{n}_{p(n,k)}\wedge t)\big)^t \big(x(s^{n}_{p(n,k+1)}\wedge t)-x(s^{n}_{p(n,k)}\wedge t)\big) \]
\[ + \sum_{k=0}^{N(\pi^n)-1} \big(y(s^{n}_{p(n,k+1)}\wedge t)-y(s^{n}_{p(n,k)}\wedge t)\big)^t \big(y(s^{n}_{p(n,k+1)}\wedge t)-y(s^{n}_{p(n,k)}\wedge t)\big)
\]
\[ + 2 \sum_{k=0}^{N(\pi^n)-1} \big(x(s^{n}_{p(n,k+1)}\wedge t)-x(s^{n}_{p(n,k)}\wedge t)\big)^t \big(y(s^{n}_{p(n,k+1)}\wedge t)-y(s^{n}_{p(n,k)}\wedge t)\big).
\]
Replacing $[x]_{d^n}$ in the above equality we obtain the following inequality:
\[\left|[x+y]_{\sigma^n}(T) - [x]_{\sigma^n}(T)\right|\]
\[ \leq \left| \sum_{k=0}^{N(\pi^n)-1} \big(y(s^{n}_{p(n,k+1)})-y(s^{n}_{p(n,k)})\big)^t \big(y(s^{n}_{p(n,k+1)})-y(s^{n}_{p(n,k)})\big)\right|
\]
\begin{equation}
 + 2\left| \sum_{k=0}^{N(\pi^n)-1} \big(x(s^{n}_{p(n,k+1)})-x(s^{n}_{p(n,k)})\big)^t \big(y(s^{n}_{p(n,k+1)})-y(s^{n}_{p(n,k)})\big)\right|. \label{eq23}
\end{equation}
We will now show that the first term of the above inequality goes to zero as $n\to \infty$.
\[ \left| \sum_{k=0}^{N(\pi^n)-1} \big(y(s^{n}_{p(n,k+1)})-y(s^{n}_{p(n,k)})\big)^t \big(y(s^{n}_{p(n,k+1)})-y(s^{n}_{p(n,k)})\big)\right|
\]
\begin{equation}
 \leq \left| \sum_{k=0}^{N(\pi^n)-1} \big(y(s^{n}_{p(n,k+1)})-y(s^{n}_{p(n,k)})\big)^p \big(y(s^{n}_{p(n,k+1)})-y(s^{n}_{p(n,k)})\big)^{2-p}\right| \label{dimention}
\end{equation}
\[\leq \sup_{k}\left| \big(y(s^{n}_{p(n,k+1)})-y(s^{n}_{p(n,k)})\big)^{2-p} \right| \left| \sum_{k=0}^{N(\pi^n)-1} \big(y(s^{n}_{p(n,k+1)})-y(s^{n}_{p(n,k)})\big)^p \right| \to 0.
\]
Using the H\"older inequality on the last term of the inequality (\ref{eq23}) we get:
\ba \left| \sum_{k=0}^{N(\pi^n)-1} \big(x(s^{n}_{p(n,k+1)})-x(s^{n}_{p(n,k)})\big)^t \big(y(s^{n}_{p(n,k+1)})-y(s^{n}_{p(n,k)})\big)\right| \label{eq.29}
\ea
\[\leq \left|\sum_{k=0}^{N(\pi^n)-1} \big(x(s^{n}_{p(n,k+1)})-x(s^{n}_{p(n,k)})\big)^t \big(x(s^{n}_{p(n,k+1)})-x(s^{n}_{p(n,k)})\big)\right|^{\frac{1}{2}}\]
\[\times \left|\sum_{k=0}^{N(\pi^n)-1} \big(y(s^{n}_{p(n,k+1)})-y(s^{n}_{p(n,k)})\big)^t \big(y(s^{n}_{p(n,k+1)})-y(s^{n}_{p(n,k) })\big)\right|^{\frac{1}{2}}.\]
Since $x$ has finite quadratic variation along $\sigma$ and  $[y]_\sigma=0$,  we have $[x,y]_\sigma=0$.
Hence we obtain $\left|[x+y]_{\sigma^n}(T) - [x]_{\sigma^n}(T)\right| \to 0$. Similarly we can obtain $\left|[x+y]_{\mathbb{T}^n}(T) - [x]_{\mathbb{T}^n}(T)\right| \to 0$. So, for the function $x+y$ we have for all $t\in [0,T]$: $[x+y]_{\mathbb{T}}(t) =[x+y]_{\sigma}(t) $ and hence $x+y\in R^{\beta}_{\pi}([0,T],\mathbb{R}^d)$. 
An upper bound for \eqref{eq.29} for $p=1$ is given by
\[\leq \left| \sum_{k=0}^{N(\pi^n)-1} \big(y(s^{n}_{p(n,k+1)-1})-y(s^{n}_{p(n,k)-1})\big)^t \big(y(s^{n}_{p(n,k+1)-1})-y(s^{n}_{p(n,k)-1})\big)^{1}\right| \]
\[\leq \sup_{k}\left\| \left(y(s^{n}_{p(n,k+1)-1})-y(s^{n}_{p(n,k)-1})\right) \right\|\times \left| \sum_{k=0}^{N(\pi^n)-1} \left\|y(s^{n}_{p(n,k+1)-1})-y(s^{n}_{p(n,k)-1}) \right\| \right|\to 0.
\]
$y\in C^0([0, T],\mathbb{R})$ has bounded variation, so the second term is finite and the first term goes to zero. Hence the above quantity goes to zero. The rest of the proof in the multidimensional case follows the same steps as above.\\
\textbf{Proof of 5.}
If $x\in R^\beta_\pi([0,T],\mathbb{R}^d)$,  there exists a dyadic subsequence or super-sequence $(d^n)_{n\geq 1}$ which satisfies Condition (i) and (ii) of Definition \ref{def.rough}. Thus there exists $N_0$ such that $|d^n|\leq 1$ for all $n\geq N_0$. Now if we fix $\gamma\geq \beta$ then for all $n\geq N_0$: $|d^n|^\gamma\leq |d^n|^{\beta}=O\left(|\pi^n|\right)$. So $x\in R^\beta_\pi([0,T],\mathbb{R}^d)$.



\begin{thebibliography}{10}

\bibitem{ananova2018}
{\sc A.~Ananova}, {\em Pathwise integration and functional calculus for paths
  with finite quadratic variation}, {PhD} thesis, Imperial College London,
  2018.

\bibitem{ananova2017}
{\sc A.~Ananova and R.~Cont}, {\em Pathwise integration with respect to paths
  of finite quadratic variation}, Journal de Math{\'e}matiques Pures et
  Appliqu{\'e}es, 107 (2017), pp.~737--757.

\bibitem{bertoin1987}
{\sc J.~Bertoin}, {\em Temps locaux et int\'egration stochastique pour les
  processus de {D}irichlet}, in S\'eminaire de {P}robabilit\'es, {XXI},
  vol.~1247 of Lecture Notes in Math., Springer, Berlin, 1987, pp.~191--205.

\bibitem{catellier2016}
{\sc R.~Catellier and M.~Gubinelli}, {\em Averaging along irregular curves and
  regularisation of {ODEs}}, Stochastic Processes and their Applications, 126
  (2016), pp.~2323 -- 2366.

\bibitem{chiu2018}
{\sc H.~Chiu and R.~Cont}, {\em On pathwise quadratic variation for c\`adl\`ag
  functions}, {Electronic Communications in Probability}, 23 (2018).

\bibitem{chiu2021}
{\sc H.~Chiu and R.~Cont},  {\em Causal functional
  calculus}, arXiv:1912.07951,  (2021).

\bibitem{cont2012}
{\sc R.~Cont}, {\em {Functional Ito Calculus and functional Kolmogorov
  equations}}, in Stochastic Integration by Parts and Functional Ito Calculus
  (Lecture Notes of the Barcelona Summer School in Stochastic Analysis, July
  2012), Advanced Courses in Mathematics, Birkhauser Basel, 2016, pp.~115--208.

\bibitem{CF2010}
{\sc R.~Cont and D.-A. Fourni\'{e}}, {\em Change of variable formulas for
  non-anticipative functionals on path space}, J. Funct. Anal., 259 (2010),
  pp.~1043--1072.

\bibitem{perkowski2019}
{\sc R.~Cont and N.~Perkowski}, {\em Pathwise integration and change of
  variable formulas for continuous paths with arbitrary regularity},
  Transactions of the American Mathematical Society, 6 (2019), pp.~134--138.

\bibitem{davis2014}
{\sc M.~Davis, J.~Obloj, and V.~Raval}, {\em Arbitrage bounds for prices of
  weighted variance swaps}, Mathematical Finance, 24 (2014), pp.~821--854.

\bibitem{davis2018}
{\sc M.~Davis, J.~Obłój, and P.~Siorpaes}, {\em Pathwise stochastic calculus
  with local times}, Ann. Inst. H. Poincaré Probab. Statist., 54 (2018),
  pp.~1--21.

\bibitem{delavega1974}
{\sc W.~F. de~La~Vega}, {\em On almost sure convergence of quadratic {B}rownian
  variation}, Ann. Probab., 2 (1974), pp.~551--552.

\bibitem{dm}
{\sc C.~Dellacherie and P.-A. Meyer}, {\em Probabilities and potential},
  vol.~29 of North-Holland Mathematics Studies, North-Holland Publishing Co.,
  Amsterdam, 1978.

\bibitem{dudley1973}
{\sc R.~M. Dudley}, {\em Sample functions of the gaussian process}, Ann.
  Probab., 1 (1973), pp.~66--103.

\bibitem{dudley2011}
{\sc R.~M. Dudley and R.~Norvai{\v{s}}a}, {\em Concrete functional calculus},
  Springer Monographs in Mathematics, Springer, New York, 2011.

\bibitem{follmer1981}
{\sc H.~F{\"o}llmer}, {\em Calcul d'{I}t\^o sans probabilit\'es}, in Seminar on
  {P}robability, {XV} ({U}niv. {S}trasbourg, {S}trasbourg, 1979/1980)
  ({F}rench), vol.~850 of Lecture Notes in Math., Springer, Berlin, 1981,
  pp.~143--150.

\bibitem{freedman}
{\sc D.~Freedman}, {\em Brownian Motion and Diffusion}, Springer, 1983.

\bibitem{FrizHairer}
{\sc P.~K. Friz and M.~Hairer}, {\em A course on rough paths}, Universitext,
  Springer, 2014.

\bibitem{geman1976}
{\sc D.~Geman and J.~Horowitz}, {\em Local times for real and random
  functions}, Duke Math. J., 43 (1976), pp.~809--828.

\bibitem{geman1980}
{\sc D.~Geman and J.~Horowitz}, {\em Occupation
  densities}, Ann. Probab., 8 (1980), pp.~1--67.

\bibitem{hanson1971}
{\sc D.~L. Hanson and F.~T. Wright}, {\em A bound on tail probabilities for
  quadratic forms in independent random variables}, Ann. Math. Statist., 42
  (1971), pp.~1079--1083.

\bibitem{imkeller2015}
{\sc P.~Imkeller and D.~J. Pr{\"o}mel}, {\em Existence of {L\'evy}'s area and
  pathwise integration}, {Communications on Stochastic Analysis}, 9 (2015).

\bibitem{karandikar1983}
{\sc R.~L. Karandikar}, {\em On the quadratic variation process of a continuous
  martingale}, Illinois J. Math., 27 (1983), pp.~178--181.

\bibitem{karandikar1995}
{\sc R.~L. Karandikar},  {\em On pathwise
  stochastic integration}, Stochastic Process. Appl., 57 (1995), pp.~11--18.

\bibitem{karandikar2014}
{\sc R.~L. Karandikar and B.~V. Rao}, {\em On quadratic variation of
  martingales}, Proc. Indian Acad. Sci. Math. Sci., 124 (2014), pp.~457--469.

\bibitem{kim2019}
{\sc D.~Kim}, {\em Local time for continuous paths with arbitrary regularity},
  arxiv,  (2019).

\bibitem{levy1940}
{\sc P.~L\'evy}, {\em Le mouvement brownien plan}, American Journal of
  Mathematics, 62 (1940), pp.~487--550.

\bibitem{levy1965}
{\sc P.~L{\'e}vy}, {\em Processus stochastiques et mouvement brownien},
  Gauthier-Villars \& Cie, Paris, 1948.

\bibitem{perkowski2015}
{\sc N.~Perkowski and D.~J. Pr{\"o}mel}, {\em Local times for typical price
  paths and pathwise {T}anaka formulas}, Electron. J. Probab., 20 (2015),
  pp.~no. 46, 15.

\bibitem{protter}
{\sc P.~E. Protter}, {\em Stochastic integration and differential equations},
  Springer-Verlag, Berlin, 2005.
\newblock Second edition.

\bibitem{schied2016}
{\sc A.~Schied}, {\em On a class of generalized {T}akagi functions with linear
  pathwise quadratic variation}, J. Math. Anal. Appl., 433 (2016),
  pp.~974--990.

\bibitem{taylor1972}
{\sc S.~J. Taylor}, {\em Exact asymptotic estimates of {Brownian} path
  variation}, Duke Math. J., 39 (1972), pp.~219--241.

\bibitem{wuermli1980}
{\sc M.~Wuermli}, {\em Lokalzeiten f\"ur martingale}, diploma thesis,
  Universit\"at Bonn, 1980.

\end{thebibliography}


\end{document}